\documentclass[10pt]{amsart}
\usepackage[dvips]{graphicx}

\newtheorem{theorem}{Theorem}[section]
\newtheorem{lemma}[theorem]{Lemma}           
\newtheorem{cor}[theorem]{Corollary}

\theoremstyle{definition}
\newtheorem{definition}[theorem]{Definition}

\newtheorem{remark}[theorem]{Remark}

\numberwithin{equation}{section}

\begin{document}
\baselineskip 14pt

\author{Hiroshi ISOZAKI}
\address{Institute of Mathematics \\
University of Tsukuba,
Tsukuba, 305-8571, Japan}
\author{Yaroslav Kurylev}
\address{Department of Mathematics \\
University College London, United Kingdom}
\author{Matti Lassas}
\address{Department of Mathematics and Statistics\\
University of Helsinki, Finland}

\title[ inverse problem on asymptotically hyperbolic orbifolds]{Spectral theory and inverse problem on asymptotically hyperbolic orbifolds}
\date{25, November, 2013}

\maketitle

\begin{abstract}
We consider an inverse problem associated with $n$-dimensional asymptotically 
hyperbolic  orbifolds $(n \geq 2)$ having a finite number of cusps and regular ends. 
By observing solutions of the Helmholtz equation at the cusp, we introduce a generalized 
$S$-matrix, and then show that it determines the manifolds with its Riemannian metric 
and the orbifold 
structure.
\end{abstract}

\begin{center}
{\large Contents}
\end{center}

\S1 Introduction

\S2 Classification of 2-dimensional hyperbolic surfacess

\S3 3-dimensional hyperbolic orbifolds

\S4 Spectral properties of the model space 

\S5 Laplace-Beltrami operators on orbifolds

\S6 Generalized S-matrix

\S7 Orbifold isomorphism


\section{Introduction}


\subsection{Assumptions on the orbifold}
 We consider an $n$-dimensional $(n \geq 2)$ connected Riemannian orbifold ${\mathcal M}$, sometimes called orbifold for the sake of simplicity, which is written as a union of open sets:
\begin{equation}
{\mathcal M} = {\mathcal K}\cup{\mathcal M}_1\cup\cdots\cup{\mathcal M}_{N+N'},
\label{S1MKM1MN}
\end{equation}
where $N \geq 1, N' \geq 0$ are integers and $\mathcal M_i\cap \mathcal M_j = \emptyset$ 
if $i \neq j$.
A part $\mathcal M_i$ is henceforth called {\it an end}.
We impose the following assumptions. 

\medskip
\noindent
{\bf (A-1)} $\ $ $\mathcal K$ {\it  is a relatively compact  $n$-dimensional orbifold}.

\medskip
\noindent
{\bf (A-2)}  \  {\it Letting $\approx$ stand for "diffeomorphic", we have
$$
\mathcal M_j \approx \left\{
\begin{split}
&M_j\times (1,\infty), \quad { for} \quad 1 \leq j \leq N, \\
& M_j\times(0,1), \quad {for} \quad N+1 \leq j \leq N+N',
\end{split}
\right.
$$
where in both cases, $M_j$  is a compact  $(n-1)$-dimensional Riemannian orbifold  whose metric is denoted by $ h_j(x,dx)$.}

\medskip
\noindent
{\bf (A-3)}
{\it If $x \in M_j^{sing}$ (see (\ref{S2Mreg}) for "$M^{sing}$"),
it has a neighbourhood $U$ such that for $1 \leq j \leq N$, $U \times (1,\infty)$ has a uniformizing cover, and for $N+1\leq j\leq N+N'$, $U\times (0,1)$ has a uniformizing cover. Namely,
there is an open neighbourhood of $\,0 \in R^{(n-1)},\, \tilde U \subset R^{(n-1)}$, a discrete
subgroup $\Gamma \subset SO(n-1)$, a surjection $\pi : \tilde U \times(1,\infty) \to U\times(1,\infty)$ or $\pi : \tilde U \times(0,1) \to U\times(0,1)$  and a metric 
$\tilde g$ such that $\pi^{\ast}g = \tilde g$ and
\begin{equation}
d\tilde s^2 = \tilde g(\tilde x, y, d\tilde x, dy)=
\frac{(dy)^2 + \tilde h_j(\tilde x,y,d\tilde x,dy)}{y^2},
\label{S1ds2regular}
\end{equation}
\begin{equation}
\tilde x \in \tilde U, \quad 
y \in \left\{
\begin{split}
& (1,\infty), \quad 1 \leq j \leq N,\\
& (0,1), \quad N+1\leq j\leq N+N',
\end{split}
\right.
\nonumber
\end{equation}
where $\tilde h_j(\tilde x,y,d\tilde x,dy)$ is a symmetric covariant tensor having 
the form : 

\noindent
(i) For $1 \leq j \leq N$, 
\begin{equation}
\widetilde h_j(\widetilde x,y,d\widetilde x, dy) = \widetilde h_j(\widetilde x,d\widetilde x).
\label{S1ds2cuspexpand}
\end{equation}
(ii) For $N+1 \leq j \leq N+N'$,
\begin{equation}
\begin{split}
& \tilde h_j(\tilde x,y,d \tilde x, dy) \\
& = \tilde h_j(\tilde x,d \tilde x)  + \sum_{p,q=1}^{n-1}a_{j,pq}(\tilde x,y)d\tilde x_p d\tilde x_q + 
\sum_{p=1}^{n-1}b_{j,p}(\tilde x,y)d\tilde x_pdy + c_j(\tilde x,y)dy^2,
\end{split}
\label{S1ds2regularexpand}
\end{equation}
where $a_{j,pq}(\tilde x,y), b_{j,p}(\tilde x,y), c_j(\tilde x,y)$ have the estimate 
for any $\alpha, \beta$
\begin{equation}
|(\partial_{\tilde x})^{\alpha}(y\partial_y)^{\beta}a(\tilde x,y)| 
\leq C_{\alpha\beta}
(1 + |\log y|)^{-\min(|\alpha|+\beta,1)-1-\epsilon_0},
\label{S1RegularAssumption}
\end{equation}
for some constant $\epsilon_0>0$. Moreover, 
$$
\gamma^* \tilde g=\tilde g, \quad {\rm for} \quad \gamma \in \Gamma,
$$
and 
$$
U \times (1,\infty) \equiv \left(\tilde U \times (1,\infty))  \right)/ \,\Gamma, \quad 
U \times (0, 1) \equiv \left(\tilde U \times (0, 1))  \right)/ \,\Gamma,
$$
where $\tilde U \times (1,\infty)$ and  $\tilde U \times (0, 1)$ are equipped with the metric $\tilde g$ and $\equiv$ stands
for the isometry.}

\medskip

A precise definition of an orbifold will be given in Subsection 1.4. 
We should note here that $C^{\infty}$ manifolds are orbifolds
(without singular points), and
that the latter may have singularities in their differential structures. 
The above assumptions allow important classes of arithmetic surface 
$\Gamma\backslash{\bf H}^n$, where ${\bf H}^n$ is the $n$-dimensional hyperbolic 
space and $\Gamma$ is a discrete subgroup of isometries on ${\bf H}^n$. If $n=2$,
all geometrically finite Fuchsian groups satisfy the assumption (A-1), (A-2), (A-3), 
where each $M_j$ is a $C^{\infty}$ manifold, as a matter of fact, $S^1$. We shall 
summerize the related results on the classification of $2-$dimensional hyperbolic 
surfaces in \S2. 
We will give a 3-dimensional example in \S3 by taking $\Gamma = SL(2,{\bf Z} + i{\bf Z})$, 
in which case $N' = 0$, and all $M_j$'s are not $C^{\infty}$ manifolds but orbifolds 
with singularities.

The orbifold at infinity, $M_j$, is called {\it a cusp} for $1 \leq j \leq N$ 
and {\it a regular infinity} for $N+1 \leq j \leq N + N'$. Sometimes $\mathcal M_j$ is 
called a cusp if $M_j$ is a cusp. We also call $\mathcal M_j$ {\it regular end} if 
$M_j$ is a regular infinity.
Let $H = - \Delta_g - (n-1)^2/4$, where $\Delta_g$ is the Laplace-Beltrami operator 
for $\mathcal M$. 
As will be explained later, it has continuous spectrum $\sigma_c(H) = [0,\infty)$, 
and the discrete spectrum $\sigma_d(H) \subset (-\infty,0)$. If $N' = 0$, $H$ may 
have embedded eigenvalues in $(0,\infty)$, which are discrete with possible 
accumulation points $0$ and $\infty$. 


\subsection{Inverse scattering from regular ends}
Let us recall known results for the case when $\mathcal M$ is a smooth Riemannian manifold 
(not orbifold) all of whose ends satisfy the assumption (A-3), i.e. are regular ends.
One can then introduce the {\it S-matrix} by observing the behavior of solutions to 
the time-dependent Schr{\"o}dinger equation or the wave equation on $\mathcal M$. An 
equivalent way is to observe the asymptotic expansion of  solutions to the Helmholtz 
equation on $\mathcal M$ in the function space $\mathcal B^{\ast}$, to be explained in \S 4.  Roughly,  $u \in \mathcal B^{\ast}$ means that $u$ behaves like 
$O(y^{(n-1)/2})$ on each end. 
One can then talk about the inverse problem. Suppose we are given two such manifolds 
$\mathcal M^{(1)}$, $\mathcal M^{(2)}$, and assume  $\mathcal M_1^{(i)}$ are regular 
ends for $i = 1,2$. Assume also that we are given, for all $k >0$, the component 
$\widehat S_{11}(k)$ of the 
$S$-matrix, which describes 
the wave coming in through $\mathcal M_1$ and going out of $\mathcal M_1$.  
Suppose that, for $\mathcal M^{(1)}$ and $\mathcal M^{(2)}$, the associated $S$-matrices 
coincide, namely, $\widehat S_{11}^{(1)}(k) = \widehat S_{11}^{(2)}(k)$ for all $k > 0$. 
If, furthermore, 
two ends $\mathcal M_1^{(1)}$ and $\mathcal M_1^{(2)}$ are known to be isometric, these 
two manifolds $\mathcal M^{(1)}$ and $\mathcal M^{(2)}$ are shown to be isometric 
(see \cite{IsKu09}).

Sa Barreto \cite{SaBa05} proved that,
 under the framework of scattering 
theory due to Melrose \cite{Me95}, two such manifolds are isometric, if the whole scattering matrix 
for all energies coincide, without assuming that one end is known to be isometric. In 
Melrose's theory of scattering metric, $(ds)^2$ is assumed to have an asymptotic expansion 
with respect to $y$ around $y = 0$, which implies that the perturbation decays 
exponentially with respect to the hyperbolic metric $\log y$.

The described results for a manifold with regular ends can be extended to the case of a manifold possessing both regular ends and  cusps. In fact, 
under the assumptions (A-1), (A-2), (A-3), where 
$M_i$ is a smooth compact Riemannian manifold, then we can get the same conclusion as above. More precisely, instead of (\ref{S1ds2cuspexpand}), for $1 \leq j \leq N$, we assume $\widetilde h_j(\widetilde x,y,d\widetilde x, dy)$ also has the form (\ref{S1ds2regularexpand}), where the coefficients of the perturbation have the estimate 
\begin{equation}
|(\widetilde D_{\tilde x})^{\alpha}(y\partial_y)^{\beta}a(\tilde x,y)| 
\leq C_{\alpha\beta}(1 + |\log y|)^{-\min(|\alpha|+\beta,1)-1-\epsilon_0}, 
\nonumber
\end{equation}
for $1 \leq j\leq N+N'$, where $\widetilde D_x = \widetilde y\partial_{\widetilde x}$ with $\widetilde y(y) = y$ for $1\leq j \leq N$, $\widetilde y(y)=1$ for $N+1\leq j \leq N+N'$.

Namely, let $\mathcal M^{(i)},\, i=1,2$,
be Riemannian manifolds of the described type and
 one of the components of their $S$-matrices associated with a regular end, say $\mathcal M^{(1)}_{N+1}$ and $\mathcal M^{(2)}_{N+1}$ coincide. Then this,  together 
with an isometry of $\mathcal M^{(1)}_{N+1}$ and $\mathcal M^{(2)}_{N+1}$,
renders the isometry of the manifolds $\mathcal M^{(i)},\, i=1, 2$.


\subsection{Main result}
The problem we address in this paper is the case in which we observe the waves coming in and going out from a cusp. 
Recall that the end $\mathcal M_1$ has a cusp at infinity.
Since the contiunuos spectrum due to the cusp is 1-dimensional, the associated S-matrix component $\widehat S_{11}(k)$ is a complex number, and it does not have enough information to determine the whole orbifold. Therefore we generalize the notion of S-matrix in the following way.

The Helmholtz equation has the following form in the cusp $\mathcal M_1$:
\begin{equation}
\left[- y^2(\partial_y^2 + \Delta_h) + (n-2)y\partial_y - \frac{(n-1)^2}{4}\right]u = k^2 u,
\label{S1equationcusp}
\end{equation}
where $k > 0$ and $\Delta_h$ is the Laplace-Beltrami operator for $M_1$. Let $\lambda_1 \leq \lambda_2 \leq \cdots$ be the eigenvalues of $- \Delta_h$ and $\varphi_1, \varphi_2, \cdots$ the associated orthnormal eigenvectors for $- \Delta_h$. Then we see that all the solution of (\ref{S1equationcusp}) has the asymptotic expansion 
\begin{equation}
u \simeq \sum_m a_my^{(n-1)/2}\varphi_m(x)e^{\sqrt{\lambda_m}y} + 
\sum_m b_m y^{(n-1)/2}\varphi_m(x)e^{-\sqrt{\lambda_m}y}
\end{equation}
as $y \to \infty$. We propose to call the operator
\begin{equation}
\mathcal S_{11}(k) : \{a_n\} \to \{b_n\}
\label{S1GneralizedSmatrix}
\end{equation}
{\it generalized S-matrix}, actually its (11) component  (see \S 6 for the precise definition). We shall show that this generalized S-matrix determines the whole orbifold $\mathcal M$. 
 
Our main result is the following.


\begin{theorem}
Suppose  two Riemannian orbifolds $\mathcal M^{(1)}$ and $\mathcal M^{(2)}$ satify the assumptions (A-1), (A-2),  (A-3). Assume that the (11) components of the generalzied scattering matrix coincide :
$$
\mathcal S_{11}^{(1)}(k) = \mathcal S_{11}^{(2)}(k), \quad \forall k > 0, \quad k^2 \not\in \sigma_p(H^{(1)})\cup\sigma_p(H^{(2)}).
$$ 
Assume also that the ends $\mathcal M_1^{(1)}$ and $\mathcal M_1^{(2)}$ are isometric. Then  $\mathcal M^{(1)}$ and $\mathcal M^{(2)}$ are isometric orbifolds.
\end{theorem}

 The main theorem 1.1 roughly means that when we send exponentially growing waves from a cusp and look at exponentially decaying waves in the same end, and if the resulting observations are the same for two asymptotically hyperbolic orbifolds, then these two orbifolds are isometric. Note that by the result \cite{Ze}, the physical scattering matrix does not determine the isometry class of the orbifold. 

Let us also remark that the numbers of ends in Theorem 1.1 are not assumed to be equal a priori. What we require is that $\mathcal M^{(1)}$ and $\mathcal M^{(2)}$ have one isometric end in common, and that the generalized scattering matrices coincide on that end.

Let us add one remark on the assumption (\ref{S1ds2cuspexpand}).
For the forward problem, i.e. the limiting absorption principle, spectral representation, the asymptotic expansion of solutions to the Helmholtz equation in $\mathcal B^{\ast}$, we have only to assume for the cusp ends ($1 \leq j \leq N$) that the coefficients of the perturbation of the metric decays like
\begin{equation}
|(y\partial_{\tilde x})^{\alpha}(y\partial_y)^{\beta}\, a(\tilde x,y)|
\leq C_{\alpha\beta}(1 + |\log y|)^{-\min (|\alpha|+\beta,1)-1-\epsilon_0},
\quad y > 1.
\label{S1forwardassumpcusp}
\end{equation}
This assumption is also sufficient for the inverse scattering from regular ends. For the inverse scattering from cusp, we have only to assume (\ref{S1ds2cuspexpand}) for the cusp end $\mathcal M_1$.

Note  that for the case of cylindrical ends, see \cite{IKL09}, the physical scattering matrix does not determine by itself the underlying manifold structure and it was necessary to introduce generalized scattering matrix. However, in that case, this generalized scattering matrix was completely determined by the physical scattering matrix.

We also note that for the potential scattering, the fact that the potential is uniquely determined by the scattering matrix in the cusp is shown in \cite{Iso04}.


\subsection{Riemannian orbifolds} \label{R-orbifolds}
Here we recall basic notions of orbifolds, see \cite{Thurs, Sa57} for further details.
A complete metric space $(\mathcal M,d)$ is called a {\it Riemannian orbifold} of dimension $n$ if, for any 
$p \in \mathcal M$,   
\begin{itemize}
\item [(i)] There exists a
radius $r(p)>0$. 
\item  [(ii)] In the 
 Riemannian ball $B_r(\widetilde p\, ;\widetilde g_p) \subset {\bf R}^n$ centered at $\widetilde p \in {\bf R}^n$ with radius $r=r(p)$ endowed with a Riemannian metric $\widetilde g_p$
 there exists a finite group $G_p \subset SO(n)$ acting by isometries
on $B_r(\widetilde p\,;\widetilde g_p)$,
\begin{equation}
\gamma^{\ast}\widetilde g_p = \widetilde g_p, \quad \forall \gamma \in G_p. 
\label{S1invariantmetric}
\end{equation}
\item  [(iii)] For all $\gamma \in G_p$ it holds that
$\gamma\cdot\widetilde p = \widetilde p$ and the action of $G_p$
on $B_r(\widetilde p\, ;\widetilde g_p)$ is faithful. 
\item  [(iv)] For the metric ball $U_p=B_r(p)\subset \mathcal M$
there exists a continuous surjection $\pi_p : B_r(\widetilde p\,;\,\widetilde g_p) \to B_r(p)$ such that 
$\pi_p(\widetilde p) = p$, moreover  for $\tilde x\in B_r(\widetilde p\,;\,\widetilde g_p)$ and
$x=\pi_p(\tilde x)$ 
it holds  $\pi_p^{-1}(x)=\{\gamma\cdot \tilde x\, ;\ \gamma\in G_p\}$.

\item  [(v)]  For  all $x,y\in B_r(p)$ the distance $d(x,y)$ satisfies
 $d(x,y)=\min\{\tilde d(\tilde x,\tilde y)\, ;\ \tilde x\in \pi_p^{-1}(x),\ \tilde  y\in \pi_p^{-1}(y)\}$, where $\tilde d$ is the Riemannian distance function on 
$B_r(\widetilde p\,;\,\widetilde g_p)$.
%
%
\end{itemize}

We say that the Riemannian ball $B_r(\widetilde p\, ;\widetilde g_p) \subset {\bf R}^n$  with a Riemannian metric $\widetilde g_p$ is the uniformizing cover of $U_p$.

\medskip
 We denote
\begin{equation}
\mathcal M_{sing} = \{p \in \mathcal M\, ;\, G_p \neq \{1\}\}, \quad
\mathcal M_{reg} = \mathcal M\setminus\mathcal M_{sing}.
\label{S2Mreg}
\end{equation}
Then,

\smallskip
\noindent
{\bf (RO-1) } {\it  $\mathcal M_{sing}$ is a closed subset of $\mathcal M$.}

\smallskip
\noindent
{\bf (RO-2)} {\it  For each $p \in \mathcal M,\, \gamma \in G_p$, $\pi_p\circ \gamma = \pi_p \ $
on $B_r(\widetilde p\, ; \widetilde g_p)$.}

\smallskip
\noindent
{\bf (RO-3)}  {\it  If $p \in M_{reg}$, $\pi_p$ is a homeomorphism from $B_r(\widetilde p\,;\,\widetilde g_p)$ to 
$U_p$.}

\noindent {\bf (RO-4) } {\it If $U_p\cap U_q \neq \emptyset$ for $p, q \in \mathcal M$, for any 
$\widetilde x \in B_r(\widetilde p\,;\,\widetilde g_p)$ and 
$\widetilde y \in B_r(\widetilde q\,;\,\widetilde g_q)$ such that 
$\pi_p(\widetilde x) = \pi_q(\widetilde y)$, there are neighborhoods 
$\widetilde x \in  V_{\widetilde x} \subset B_r(\widetilde p\,;\,\widetilde g_p)$, 
$\widetilde y \in  V_{\widetilde y} \subset B_r(\widetilde q\,;\,\widetilde g_q)$ and a 
Riemannian isometry $\psi :  V_{\widetilde x} \to  V_{\widetilde y}$ such that 
$\psi(\widetilde x) = \widetilde y$.}

\medskip
The elements in
$\mathcal M_{sing}$ are called {\it singular points} of $\mathcal M$.

\medskip
Two Riemannian orbifolds $\mathcal M$, $\mathcal M'$ are said to be {\it isometric} if they satisfy the following conditions:

\smallskip
\noindent
{\bf (I-1)} {\it There exists a homeomorphism 
$f : \mathcal M \to \mathcal M'$ such that $f(\mathcal M_{reg}) = \mathcal M'_{reg}$, $f(\mathcal M_{sing}) = \mathcal M'_{sing},$ and 
$f\big|_{\mathcal M_{reg}} : \mathcal M_{reg} \to \mathcal M'_{reg}$ is a Riemannian isometry. }\\
\smallskip
\noindent
{\bf (I-2)} {\it  For any $p \in \mathcal M_{sing}$
and $p'=f(p)\in \mathcal M'_{sing}$ 
there exist $r>0$ and  a uniformizing cover
$B_r(\widetilde p\,;\,\widetilde g_p)$ of
$B_r(p)\subset \mathcal M$ and 
 a uniformizing cover
$B_r(\widetilde {f(p)}\,;\,\widetilde g_p')$ of
$B'_r({f(p)})\subset \mathcal M'$,
and  a Riemannian isometry 
$\widetilde F_p : B_r(\widetilde p\,;\,\widetilde g_p) \to B_r(\widetilde {f(p)}\,;\,
\widetilde{g'_{f(p)}})$ such that 
\begin{equation}
f\circ \pi_p = \pi'_{f(p)}\circ \widetilde F_p,\quad \hbox{on }\quad
B_r(\widetilde p\,;\,\widetilde g_p).
\label{S1fcircpip}
\end{equation}
 Moreover, there exists a group isomorphism 
$i_p : G_p \to G'_{f(p)}$ between the group $G_p$ acting on  $B_r(\widetilde p\,;\,\widetilde g_p)$ and the group $G'_{f(p)}$ acting on  $B_r(\widetilde {f(p)}\,;\,
\widetilde{g'_{f(p)}})$ such that 
for all $\gamma \in G_p$ we have
\begin{equation}
\widetilde F_p\circ\gamma = (i_p\gamma)\circ\widetilde F_p
\quad \hbox{on }
B_r(\widetilde p\,;\,\widetilde g_p).
\label{S1Fpgamma}
\end{equation}
}

\medskip
To introduce a function calculus on Riemannian orbifolds, we use a uniformising cover.
We put $\widetilde U_p = B_r(\widetilde p\,;\,\widetilde g_p)$ for brevity.
First let us note that there is a 1 to 1 correspondence between a (one-valued) function $f$ on $U_p$ and a $G_p$-invariant function $\widetilde f$ on $\widetilde U_p$ : 
\begin{equation}
f(\pi_p(z)) = \widetilde f(z), \quad z \in \widetilde U_p ; \quad f(x) = \widetilde f(\pi_p^{-1}(x)), \quad x \in U_p.
\label{S4fandwidetildef}
\end{equation}
Note also that a function $\widetilde f$ on $\widetilde U_p$ is $G_p$-invariant if and only if there exists a function $F$ defined on $\widetilde U_p$ such that
$$
\widetilde f(z) = \frac{1}{{\#} (G_p)}\sum_{\gamma\in G_p}F(\gamma\cdot z), \quad z \in \widetilde U_p.
$$
We say that $f \in C^{\infty}(U_p)$ if $\widetilde f$ in (\ref{S4fandwidetildef}) satisfies
$\widetilde f \in C^{\infty}(\widetilde U_p)$, and $f \in H^m(U_p)$ = the Sobolev space of order $m$, if  $\widetilde f \in H^m(\widetilde U_p)$. 
Then we have : $f \in C^{\infty}_0(U_p)$ if and only if $ \widetilde f \in C^{\infty}_0(\widetilde U_p)$. 

The integral of a function $f$ over $U_p$ is defined by
\begin{equation} 
\int_{U_p}fdU_p = \frac{1}{\# (G_p)}\int_{\widetilde U_p}\widetilde fd\widetilde U_p,
\label{S4intoverUp}
\end{equation}
where $d\widetilde U_p$ is the Riemannian volume element of $\widetilde U_p$. 

If $\mathcal M$ is compact, the set of singular points is also compact. We can then construct an open covering of $\mathcal M$, $\cup_{i=1}^mV_i = \mathcal M$, such that each $V_i$ is one of the above $U_p$. By using the partition of unity $\{\chi_i\}_{i=1}^m$ subordinate to this covering, we define the integral over $\mathcal M$ by
$$
\int_{\mathcal M}fd\mathcal M = \sum_{i=1}^m\int_{V_i}\chi_ifdV_i.
$$
When $\mathcal M$ is non-compact, the situation may be more complicated. We restrict ourselves to 
the case in which the assumptions in \S1 are satisfied. Moreover, for the sake of simplicity of 
notations, we consider the case where there is only one end $\mathcal M_1$ in the assumption (A-2). 
Recall that $\mathcal M_1 = M_1\times(1,\infty)$, $M_1$ being a compact orbifold of dimension $(n-1)$, 
hence the singular points of $M_1$ form a compact set. Then as above one can construct a 
covering $\{V_j\}_{j=1}^{m'}$ of $M_1$ and a partition of unity $\{\psi_j\}_{j=1}^{m'}$ of $M_1$ 
such that 
$\psi_j \in C_0^{\infty}(V_j)$. We take $\chi(y) \in C^{\infty}(0,\infty)$ such that $\psi(y) = 0$ 
for $y < 3/2$, $\chi(y) = 1$ for $y > 2$, and put 
$\varphi_j(x,y) = \psi_j(x)\chi(y)$, $x \in M_1$. We next take into account of the compact part 
$\mathcal K\cup(\mathcal M_1\cap\{y \leq 2\})$ and add another system of functions $\tilde\varphi_{1},\cdots,\tilde\varphi_{m''}$, so that
$\tilde\varphi_{1},\cdots,\tilde\varphi_{m''}$ together with 
$\varphi_{1},\cdots,\varphi_{m'}$ form a partition of 
unity of $\mathcal M$. We rearrange them and denote by $\chi_1,\chi_2,\cdots, \chi_{m},\, m=m'+m''$.
The inner product of $L^2(\mathcal M)$ is then defined by
$$
(f,g) = \sum_{j=1}^m\int_{U_j}\chi_jf\overline{g}\, dU_j,
$$
Rellich's selection theorem and Sobolev's imbedding theorem are extended to orbifolds 
(see \cite{Chian90}).

By shrinking $\widetilde U_p$, we can assume that each $\widetilde U_p$ is in  one coordinate patch. Let $z = (z_1,\cdots,z_n)$ be the associated local coordinate on $\widetilde U_p$ and $\sum_{i,j=1}^n\widetilde g_{ij}(z)dz_idz_j$  the $G_p$-invariant Riemannian metric on $\widetilde U_p$.
Then, letting $x = \pi_p(z)$, the Laplace-Beltrami operator $\Delta_g$ for the orbifold $\mathcal M$ is written as
$$
\Delta_g f(x) = \left(
\frac{1}{\sqrt{\widetilde g}}\sum_{i,j=1}^n\frac{\partial}{\partial z_i}\sqrt{\widetilde g}\widetilde{g}^{ij}\frac{\partial}{\partial z_j}\widetilde f\right)\left(\pi^{-1}_p(x)\right), \quad 
\widetilde g = \det\left(\widetilde{g}_{ij}\right).
$$
As in the case of Riemannian manifolds, $\Delta_g$ is symmetric on $C_0^{\infty}(\mathcal M)$. If $\mathcal M$ is compact, one can show that $-\Delta_g$ with domain $H^2(\mathcal M)$ is self-adjoint with discrete spectrum $0 = \mu_1 \leq \mu_2 \leq \cdots \to \infty$ (for some results regarding the spectral theory of the Laplacian on compact orbifolds, see \cite{DGGW08} and \cite{Far01}).


\subsection{Plan of the paper}

In \S2, we recall a classical classification theorem of 2-dimensional hyperbolic spaces, and in \S3 a basic example of the 3-dimensional orbifold. General properties of the Laplace-Beltrami operator on orbifold are discussed in \S 4. In \S 5 and \S6, we study spectral properties of the Laplace-Beltrami operator of our orbifolds. Generalized S-matrix is defined in \S 7. We shall prove Theorem 1.1 in \S 8.

In our lecture notes \cite{IsKu09}, we have studied
spectral properties of asymptotically hyperbolic manifolds in detail. 
The case of orbifold requires essentially no change, since by the spectral decomposition of the Laplace-Beltrami operator $\Delta_h$ of the associated orbifold at infinity, we are reduced to the 1-dimensional problem, which is just the case we have studied in \cite{IsKu09}. However, for the reader's convenience, we shall develop here the spectral theory in as much detail as possible, so that 
this paper becomes self-contained. In order not to make this paper too long, we have omitted routine computational parts, which can be done by straightforward calculation, or seen in \cite{IsKu09}.

The notation used in this paper is standard. For a manifold (or orbifold) $\mathcal M$ with the volume element $dV$, and a Banach space $X$, $L^2(\mathcal M;X;dV)$ denotes the Hilbert space of $X$-valued $L^2$-functions with respect to the measure $dV$. We put $L^2(\mathcal M;dV) = L^2(\mathcal M;{\bf C};dV)$.
$H^m(\mathcal M)$ denotes the usual Sobolev space of order $m$ on $\mathcal M$.
 For two Banach spaces $X$ and $Y$, ${\bf B}(X;Y)$ is the set of all bounded linear operators from $X$ to $Y$. 

At the end of this introduction, we would like to note that in the very recent years spectral, scattering and inverse scattering for asymptotically hyperbolic manifolds and manifolds with singularities have attracted much interest in the mathematical community, see e.g. 
\cite{BP11}, \cite{GuMa}, \cite{GuNau}, \cite{GuSaBa08},  \cite{MaVa}, \cite{Perry1}, \cite{Vasy}, \cite{Vasy2}.


\section{Classification of 2-dimensional hyperbolic surfaces}
A hyperbolic manifold $\mathcal M$ is a complete Riemannian manifold whose sectional curvatures 
are $-1$, and it is identified with $\Gamma\backslash{\bf H}^n$, where ${\bf H}^n$ is the 
hyperbolic space and $\Gamma$ is a discrete group of isometries on ${\bf H}^n$ which may have 
fixed points. In this case, actually, $\Gamma\backslash{\bf H}^n$  is a Riemannian orbifold. 
In this section, we recall the basic facts about $\Gamma\backslash{\bf H}^n$ in the 2-dimensional case.  


\subsection{Fuchsian groups}
The upper-half space model of 2-dimensional hyperbolic space ${\bf H}^2$ is  ${\bf C}_+ = \{z = x + iy \, ; y > 0\}$ equipped with the metric
\begin{equation}
ds^2 = \frac{(dx)^2 + (dy)^2}{y^2}.
\label{S2Metric}
\end{equation}
The infinity of ${\bf H}^2$ is 
$$
\partial{\bf C_+} = \{(x,0) \, ; x \in {\bf R}\}\cup{\infty} = 
{\bf R}\cup{\infty}.
$$
It admits an action of $SL(2,{\bf R})$ defined by
\begin{equation}
SL(2,{\bf R})\times{\bf C}_+ \ni (\gamma,z) \to \gamma\cdot z:= \frac{az + b}{cz + d}, \quad \gamma = \left(
\begin{array}{cc}
a & b \\
c & d
\end{array}
\right),
\label{S2Action}
\end{equation}
where the right-hand side, {\it  M{\"o}bius transformation}, is an isometry on ${\bf H}^2$. The mapping : $\gamma \to \gamma\cdot$ is 2 to 1, and the group of M{\"o}bius transformations is isomorphic to 
$PSL(2,{\bf R}) = SL(2,{\bf R})/\{\pm I\}$.
This action is classified into 3 categories :
\begin{equation}
\begin{split}
elliptic &\Longleftrightarrow {\rm there \ is \ only \ one \ fixed \ point \ in} \ {\bf C}_+ \\
&\Longleftrightarrow |{\rm tr}\, \gamma| < 2, \\
parabolic &\Longleftrightarrow {\rm there \ is \ only \ one \ degenerate \ fixed \ point \ on} \ \partial{\bf C}_+ \\
&\Longleftrightarrow |{\rm tr}\, \gamma| = 2, \\
hyperbolic &\Longleftrightarrow {\rm there \ are \ two \  fixed \ points \ on} \ \partial{\bf C}_+ \\
&\Longleftrightarrow |{\rm tr}\, \gamma| > 2.
\end{split}
\nonumber
\end{equation}

Let $\Gamma$ be a discrete subgroup, {\it Fuchsian group}, of $SL(2,{\bf R})$, and 
$\mathcal M = \Gamma\backslash{\bf H}^2$ be its fundamental domain by the action (\ref{S2Action}). 
$\Gamma$ is said to be {\it geometrically finite} if $\mathcal M$ can be chosen to be a 
finite-sided convex polygon. The sides are then the geodesics in ${\bf H}^2$. The geometric 
finiteness is equivalent to that $\Gamma$ is finitely generated. If $\Gamma$ is a Fuchsian group, 
 the quotient space $\Gamma\backslash{\bf H}^2$ is called a hyperbolic surface.
Let us give two simple 
but fundamental examples. 

\subsubsection{Parabolic cyclic group}
Consider the cyclic group $\Gamma$ generated by the action $z \to z + 1$. This is parabolic with 
the only fixed point $\infty$. The associated fundamental domain is 
$\mathcal M = [-1/2,1/2]\times(0,\infty)$,
which is a hyperbolic manifold with metric (\ref{S2Metric}). 
It has two infinities : $[-1/2,1/2]\times \{0\}$ and $\infty$. The part 
$[-1/2,1/2]\times(0,1)$ has an infinite volume. 
The part $[-1/2,1/2]\times(1,\infty)$ has a finite volume, and is called the {\it cusp}. 
The sides $x = \pm 1/2$ are geodesics. The quotient manifold is diffeomorphic to 
$S^1\times(-\infty,\infty)$.

\subsubsection{Hyperbolic cyclic group}
Another simple example is the cyclic group generated by the hyperbolic action 
$z \to \lambda z$, $\lambda > 1$. The sides of the fundamental domain 
$\mathcal M = \{1 \leq |z| \leq \lambda\}$ are semi-circles orthogonal to $\{y = 0\}$, 
which are geodesics. The quotient manifold is diffeomorphic to $S^1\times(-\infty,\infty)$. 
It can be parametrized by $(t,r)$, where $t \in {\bf R}/(\log\lambda){\bf Z}$ and $r$ is the 
signed distance from the segment $\{(0,t)\, ; 1 \leq t \leq \lambda\}$. The metric is then 
written as
\begin{equation}
ds^2 = (dr)^2 + \cosh^2r\,(dt)^2.
\label{eq:Chap3Sect1Funnelmetric}
\end{equation}
The part $r > 0$ (or $r < 0)$ is called the {\it funnel}. Letting $y = 2e^{-r},\, 0< y<1$, one 
can rewrite (\ref{eq:Chap3Sect1Funnelmetric}) as
\begin{equation}
ds^2 = \Big(\frac{dy}{y}\Big)^2 + \Big(\frac{1}{y} + \frac{y}{4}\Big)^2(dt)^2.
\nonumber
\end{equation}
Therefore, the funnel can be regarded as a perturbation of the infinite volume part 
$[-1/2,1/2]\times(0,1)$ of the fundamental domain for the parabolic cyclic group.


\subsection{Classification of hyperbolic surfaces}
The set of the limit points of a Fuchsian group $\Gamma$, denoted by $\Lambda(\Gamma)$, is 
defined as follows : $w \in \Lambda(\Gamma)$ if there exist $z_0 \in {\bf C}_+$ and 
$I \neq\gamma_n \in \Gamma$ such that $\gamma_n\cdot z_0 \to w$.
 Since $\Gamma$ acts discontinuously on ${\bf C}_+$, $\Lambda(\Gamma) \subset \partial{\bf H}^2$. 
 There are only 3 possibilities.
\begin{itemize}
\item 
({\it Elementary}) : 
$\Lambda(\Gamma)$ is a finite set.

\item
({\it The 1st kind}) :
$\ \Lambda(\Gamma) = \partial{\bf H}$.

\item
({\it The 2nd kind}) :
$\ \Lambda(\Gamma)$ is a perfect (i.e. every point is an accumulation point), nowhere dense set of $\partial{\bf H}$.
\end{itemize}
 
Any elementary group is either cyclic or is conjugate in $PSL(2,{\bf R})$ to a group generated by $\gamma\cdot z = \lambda z$, $(\lambda > 1)$, and $\gamma'\cdot z = - 1/z$ (see \cite{Ka92}, Theorem 2.4.3).

For non-elementary case, we have the following theorem (\cite{Bor07}, Theorem 2.13). 


\begin{theorem}
Let ${\mathcal M} = \Gamma\backslash{\bf H}^2$ be a non-elementary geometrically finite hyperbolic surface. Then there exists a compact subset $\mathcal K$ such that $\mathcal M \setminus {\mathcal K}$ is a finite disjoint union of cusps and funnels.
\end{theorem}

Other important theorems are the following (see \cite{Ka92}, Theorems 4.5.1, 4.5.2 and 4.1.1).


\begin{theorem}
A Fuchsian group is of the 1st kind if and only if its fundamental domain has a finite area.
\end{theorem}


\begin{theorem}
A Fuchsian group of the 1st kind is geometrically finite.
\end{theorem}

 For the Fuchsian group of the 1st kind, therefore, the ends of its fundamental domain are always cusps. 
 In this case, usually it is compactified around parabolic fixed points and made a compact 
 Riemann surface. The automorphic functions associated with this group turn out to be algebraic 
 functions on this Riemann surface.


\subsection{Orbifold structure} 
We need to be careful about the analytic structure around the elliptic fixed points. Under the 
assumption of geometric finiteness, they are finite in the fundamental domain. Let $p$ be an 
elliptic fixed point, and ${\rm Iso}(p)$  the finite cyclic isotropy group of $p$. Assume 
that ${\#} {\mathcal I}(p) = n$. Then its generator $\gamma$ satisfies
\begin{equation}
\frac{w - p}{w - \overline{p}} = e^{2\pi i/n}\frac{z - p}{z - \overline{p}}, \quad w = \gamma\cdot z.
\nonumber
\end{equation}
By the linear transformation 
\begin{equation}
T(z) = (z - p)/(z - \overline{p}),
\nonumber
\end{equation}
 it is written as
\begin{equation}
\gamma\cdot = T^{-1}\kappa T, \quad \kappa = e^{2\pi i/n}.
\nonumber
\end{equation}
Therefore, ${\rm Iso}(p)$ is isomorphic to the rotation group by the angle $2\pi/n$. We introduce 
the local coordinates $\varphi_p(\pi(z))$ around $p$ by
\begin{equation}
\zeta := \varphi_p(\iota(z)) = T(z)^n,
\nonumber
\end{equation}
where $\iota$ is the canonical projection 
\begin{equation}
\iota : {\bf H}^2 \ni z \to [z] \in \Gamma\backslash{\bf H}^2,
\nonumber
\end{equation}
with $[z] = \{g\cdot z\, ; \, g \in \Gamma\}$.
Since $n \geq 2$, the Riemannian metric (\ref{S2Metric}) has a singularity at $p$. Identifying $z$ and $\iota(z)$, we have
$$
z = \frac{p - \overline{p}\zeta^{1/n}}{1 - \zeta^{1/n}} = p + (p-{\overline p})\zeta^{1/n} + \cdots.
$$
Therefore
\begin{equation}
\frac{(dx)^2 + (dy)^2}{y^2} = \frac{dz\,d\overline{z}}{({\rm Im}\,z)^2} = 
\frac{\big|dz/d\zeta\big|^2}{({\rm Im}\,z)^2}d\zeta d\overline{\zeta},
\nonumber
\end{equation}
\begin{equation}
\left|\frac{dz}{d\zeta}\right|^2 = O(|\zeta|^{-\lambda}), \quad \lambda = 2 - \frac{2}{n}.
\nonumber
\end{equation}
Note that $1 \leq  \lambda < 2$.
The volume element and the Laplace-Beltrami operator are rewritten as
\begin{equation}
\frac{dx\wedge dy}{y^2} = \frac{i}{2y^2}dz\wedge d\overline{z} = \frac{i\big|dz/d\zeta\big|^2}{2({\rm Im}\,z)^2}d\zeta\wedge d\overline{\zeta},
\nonumber
\end{equation}
\begin{equation}
y^2\big(\partial_x^2 + \partial_y^2\big) = 4({\rm Im}\,z)^2\frac{\partial^2}{\partial z\partial\overline{z}} = 
\frac{4({\rm Im}\,z)^2}{\big|dz/d\zeta|^2}\frac{\partial^2}{\partial\zeta\partial\overline{\zeta}}.
\nonumber
\end{equation}
Both of them have singularities at $p$. However, for $f, g$, $C^{\infty}$-functions supported near $p$, we have
\begin{equation}
\int_{\Gamma\backslash{\bf H}}y^2\big(\partial_x^2 + \partial_y^2\big)f\cdot g\,\frac{dxdy}{y^2} = 2i
\int \frac{\partial^2}{\partial\zeta\partial\overline{\zeta}}f\cdot g\, d\zeta d\overline{\zeta}.
\nonumber
\end{equation}
What is important is that the singularity of the volume element and that of the Laplace-Beltrami operator cancel.


\subsection{Example} 
Let $\Gamma = SL(2,{\bf Z})$. Then 
\begin{equation}
SL(2,{\bf Z})\backslash{\bf H}^2 = \{z \in {\bf C}_+ \, ; \, |z| \geq 1, |{\rm Re}\,z| \leq 1/2\}.
\nonumber
\end{equation}
The fixed points are $i, e^{\pi i/3}, e^{2\pi i/3}$, and the isotropy groups are\begin{equation}
{\rm Iso}(i) = \Big<\left(
\begin{array}{cc}
0 & -1 \\
1 & 0
\end{array}
\right)\Big> = {\bf Z}_2,
\nonumber
\end{equation}
\begin{equation}
{\rm Iso}(e^{\pi i/3}) = \Big\langle\left(
\begin{array}{cc}
0 & -1 \\
1 & -1
\end{array}
\right)\Big\rangle = {\bf Z}_3,
\nonumber
\end{equation}
\begin{equation}
{\rm Iso}(e^{2\pi i/3}) = \Big<\left(
\begin{array}{cc}
-1 & -1 \\
1 & 0
\end{array}
\right)\Big> = {\bf Z}_3,
\nonumber
\end{equation}
where $\langle g\rangle$ denotes the cyclic group generated by $g$, and ${\bf Z}_n = {\bf Z}/n{\bf Z}$.


\section{3-dimensional hyperbolic orbifolds}

\subsection{Kleinian group}
The upper-half space model of 3-dimensional hyperbolic space ${\bf H}^3$ is  ${\bf R}_+^3 = \{(x_1,x_2,y) \, ; y > 0\}$ equipped with the metric
\begin{equation}
ds^2 = \frac{(dx_1)^2 + (dx_2)^2 + (dy)^2}{y^2}.
\nonumber
\end{equation}
The infinity of ${\bf H}^3$ is 
$$
\partial{\bf H}^3 = \{(x_1,x_2,0) \, ; (x_1,x_2) \in {\bf R}^2\}\cup{\infty} = 
{\bf R}^2\cup{\infty}.
$$
We represent points in ${\bf H}^3$ by quarternions :
\begin{equation}
(x_1,x_2,y) \longleftrightarrow x_1{\bf 1} + x_2{\bf i} + y{\bf j},
\nonumber
\end{equation}
\begin{equation}
\bf 1 = 
\left(
\begin{array}{cc}
1 & 0 \\
0 & 1
\end{array}
\right), \quad
\bf i = 
\left(
\begin{array}{cc}
i & 0 \\
0 & -i
\end{array}
\right), \quad
\bf j = 
\left(
\begin{array}{cc}
0 & 1 \\
-1 & 0
\end{array}
\right), \quad
\bf k = 
\left(
\begin{array}{cc}
0 & i \\
i & 0
\end{array}
\right),
\nonumber
\end{equation}
It is convenient to identify $x_1{\bf 1} + x_2{\bf i}$ with  $z = x_1 + ix_2 \in {\bf C}$. Then 
\begin{equation}
x_1{\bf 1} + x_2{\bf i} + y{\bf j} = 
\left(
\begin{array}{cc}
z & y \\
- y & \overline{z}
\end{array}
\right),
\nonumber
\end{equation}
which is denoted by $\zeta:= z + yj$ below.
${\bf H}^3$ admits the action of $SL(2,{\bf C})$, {\it  M{\"o}bius transformation}, defined by
\begin{equation}
SL(2,{\bf C})\times{\bf H}^3 \ni (\gamma,\zeta) \to \gamma\cdot \zeta:= (a\zeta + b)(c\zeta + d)^{-1}, \quad \gamma = \left(
\begin{array}{cc}
a & b \\
c & d
\end{array}
\right).
\nonumber
\end{equation}
Note that by the above identification,
$$
a\zeta + b = 
\left(
\begin{array}{cc}
a & 0 \\
0 & \overline{a}
\end{array}
\right)
\left(
\begin{array}{cc}
z & y \\
-y & \overline{z}
\end{array}
\right) + 
\left(
\begin{array}{cc}
b & 0 \\
0 & \overline{b}
\end{array}
\right)
 = 
\left(
\begin{array}{cc}
az + b &ay \\
- \overline{ay} & \overline{az + b}
\end{array}
\right).
$$
We have, therefore,
\begin{equation}
\gamma\cdot\zeta = \frac{1}{|cz+d|^2+ |c|^2y^2}
\left(
\begin{array}{cc}
w & y \\
-y & \overline{w}
\end{array}
\right), \quad
w = (az + b)\overline{(cz + d)} + a\overline{c}|y|^2.
\label{S3action2}
\end{equation}
This is an isometry on ${\bf H}^3$. The mapping : $\gamma \to \gamma\cdot$ is 2 to 1, and the group of M{\"o}bius transformations is isomorphic to 
$PSL(2,{\bf C}) = SL(2,{\bf C})/\{\pm I\}$.
They are classified into 4 categories :
\begin{equation}
\begin{split}
elliptic
&\Longleftrightarrow {\rm tr}\,\gamma \in {\bf R}, \quad |{\rm tr}\, \gamma| < 2, \\
parabolic 
&\Longleftrightarrow {\rm tr}\,\gamma \in {\bf R}, \quad |{\rm tr}\, \gamma| = 2, \\
hyperbolic 
&\Longleftrightarrow {\rm tr}\,\gamma \in {\bf R}, \quad |{\rm tr}\, \gamma| > 2,  \\
loxodromic &\Longleftrightarrow {\rm tr}\,\gamma \not\in {\bf R}.
\end{split}
\nonumber
\end{equation}

A subgroup $\Gamma$ of $SL(2,{\bf C})$ is called {\it Kleinian group} if $\Gamma$ is discrete 
in $SL(2,{\bf C})$. 
A Kleinian group $\Gamma$ is said to be {\it torsion free} if no element of $\Gamma\setminus\{1\}$ 
has a fixed point in ${\bf H}^3$. This is equivalent to that any element $1 \neq \gamma \in \Gamma$ 
has an infinite order. Let $\mathcal M = \Gamma\backslash{\bf H}^3$ be the fundamental domain for 
the action (\ref{S2Action}).
The following theorem is well-known, see e.g. \cite{MaTa98}, p. 28, Theorem 1.18.

\begin{theorem}
For any complete 3-dim. hyperbolic manifold $\mathcal M$, there is a torsion free Kleinian group $\Gamma$ such that $\mathcal M = \Gamma\backslash{\bf H}^3$. Conversely, for any torsion free Kleinian group $\Gamma$, $\Gamma\backslash{\bf H}^3$ is a complete 3-dim. hyperbolic manifold.
\end{theorem}

Therefore, the existence of fixed points is essential for the construction of orbifolds. 


\begin{lemma}
An element $1 \neq \gamma \in SL(2,{\bf C})$ has a fixed point in ${\bf H}^3$ if and only if $\gamma$ is elliptic. In this case, $\gamma$ has two fixed points $p, q \in \partial{\bf H}^3$, and all the points on the geodesic $C_{p,q}$ joining $p$ and $q$ are left fixed by the action of $\gamma$. Moreover, $\gamma$ is a rotation around $C_{p,q}$.
\end{lemma}

For the proof, see \cite{MaTa98}, p. 27, Proposition 1.16 and \cite{EGM97}, p. 34, Proposition 1.4.
If $\zeta = z + yj \in {\bf H}^3$ is a fixed point of $\gamma \in SL(2,{\bf C})$, we have by (\ref{S3action2})
\begin{equation}
\frac{(az + b)\overline{(cz + d)} + a\overline{c}y^2}{|cz + d|^2 + |cy|^2} = z, \quad
\frac{y}{|cz + d|^2 + |cy|^2} = y,
\label{S3fixedpoints}
\end{equation}


\subsection{Picard group}
As an example, we consider the Picard group
\begin{equation}
\Gamma = PSL(2,{\bf Z}[i]) = \left\{
\left(
\begin{array}{cc}
a & b \\
c & d
\end{array}
\right)\, ;\, a, b, c, d \in {\bf Z}[i], \ ad - bc = 1\right\},
\label{S3Picard}
\end{equation}
where ${\bf Z}[i] = {\bf Z} + i{\bf Z}$, the ring of Gaussian integers. 
The following lemma is proven in \cite{EGM97}, p. 324, Proposition 3.9.


\begin{lemma}
(1) $PSL(2,{\bf Z}[i])$ is generated by
\begin{equation}
\gamma_1 = 
\left(
\begin{array}{cc}
1 & 0 \\
1 & 1
\end{array}
\right), \quad
\gamma_2 = 
\left(
\begin{array}{cc}
0 & -1 \\
1 & 0
\end{array}
\right), \quad
\gamma_3 = 
\left(
\begin{array}{cc}
1 & 0 \\
i & 1
\end{array}
\right).
\nonumber
\end{equation}
(2) The fundamental domain of $PSL(2,{\bf Z}[i])$ is
\begin{equation}
\mathcal M = \Gamma\backslash{\bf H}^3 = 
\left\{ z + yj\, ; \,
|{\rm Re}\, z | \leq \frac{1}{2}, \ 0 \leq {\rm Im}\, z \leq \frac{1}{2}, \ 
|z|^2 + y^2 \geq 1\right\}.
\nonumber
\end{equation}
(3) The vertices of $\mathcal M$ are 
\begin{equation}
\infty, \quad -\frac{1}{2} + \frac{\sqrt3}{2}j, \quad \frac{1}{2} + \frac{\sqrt3}{2}j,\quad
\frac{1}{2} + \frac{1}{2}i + \frac{\sqrt2}{2}j, \quad
- \frac{1}{2} + \frac{1}{2}i + \frac{\sqrt2}{2}j.
\nonumber
\end{equation}
\end{lemma}

\begin{figure}[hbtp]
\centering
\includegraphics[width=9cm]{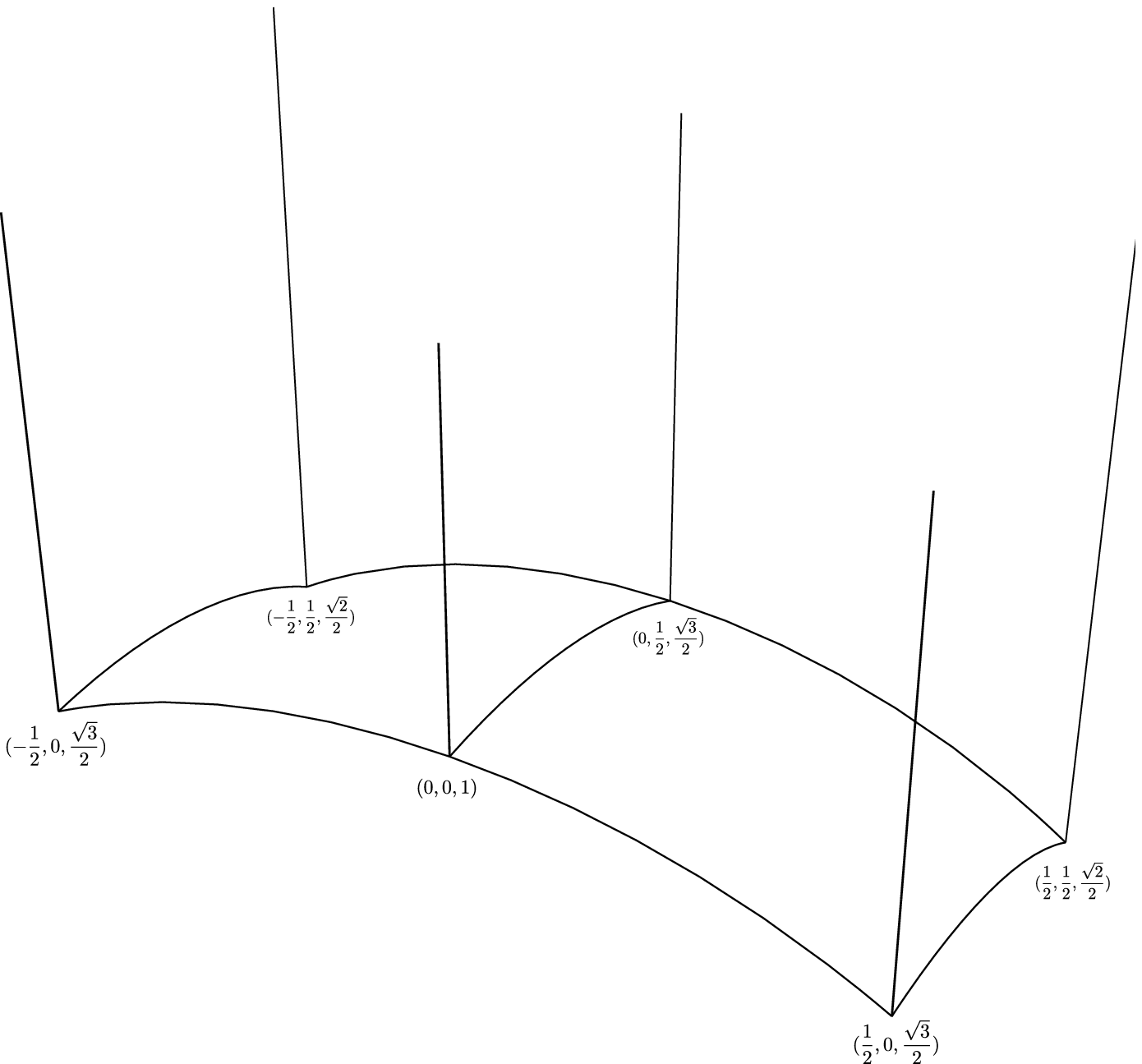}
\caption{Fundamental domain for $PSL(2,{\bf Z}[i])$ (1)}
\label{fig:math2}
\end{figure}

\begin{figure}[hbtp]
\centering
\includegraphics[width=9cm]{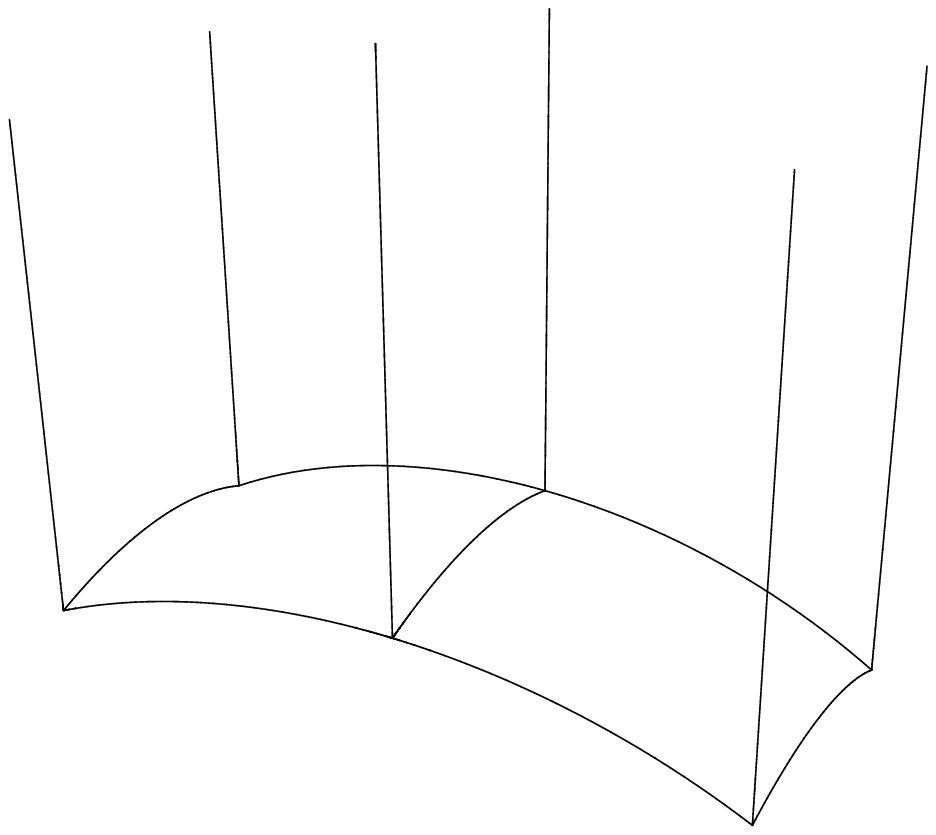}
\caption{Fundamental domain for $PSL(2,{\bf Z}[i])$ (2)}
\label{fig:math1}
\end{figure}

\medskip
Note that
\begin{equation}
g_1 = 
\left(
\begin{array}{cc}
0 & -1 \\
1 & 0
\end{array}
\right), \quad
g_2 = 
\left(
\begin{array}{cc}
i & -1 \\
0 & -i
\end{array}
\right), \quad
g_3 = 
\left(
\begin{array}{cc}
-i & 0 \\
0 & i
\end{array}
\right), \quad
h = 
\left(
\begin{array}{cc}
1 & 1 \\
0 & 1
\end{array}
\right)
\nonumber
\end{equation}
are another system of generators. In fact, we have the following relations :
$$
\gamma_1 = hg_1h, \quad \gamma_2 = g_1, \quad \gamma_3 = g_1^{-1}g_2g_3g_1.
$$
By (\ref{S3action2}), $g_1\cdot$ is the inversion with respect to the unit sphere
$$
g_1\cdot(z + yj) = \frac{-\overline{z}}{|z|^2 + y^2} + \frac{y}{|z|^2 + y^2}j,
$$
$g_2$ is a shifted horizontal inversion,
$$
g_2\cdot (z+yj)= -z-i +yj,
$$
$g_3\cdot$ is the horizontal inversion
$$
g_3\cdot(z + yj) = - z + yj,
$$
$g_3\cdot g_2\cdot$ is the translation to the $\bf i$-direction
$$
g_3\cdot g_2\cdot(z + yj) = z + i + yj,
$$
and $h\cdot$ is the translation to the $\bf 1$-direction
$$
h\cdot(z + yj) = z + 1 + yj.
$$ 
We note that the action by $g_2, g_3$ and $h$ do not affect $y$.

The boundary of the fundamental domain is split into 8 parts, $S_{j}, 1 \leq j \leq 8$ :
\begin{equation}
\begin{split}
S_{1} & = {\mathcal M}\cap \{{\rm Re}\,z = - 1/2\}, \\
S_{2} & = {\mathcal M}\cap \{{\rm Re}\,z =  1/2\}, \\
S_{3} &= {\mathcal M}\cap \{-1/2 \leq {\rm Re}\,z \leq 0, \, {\rm Im}\, z = 0\}, \\
S_{4} &= {\mathcal M}\cap \{0 \leq {\rm Re}\,z \leq 1/2, \, {\rm Im}\, z = 0\}, \\
S_{5} &= {\mathcal M}\cap \{-1/2 \leq {\rm Re}\,z \leq 0, \, {\rm Im}\, z = 1/2\}, \\
S_{6} &= {\mathcal M}\cap \{0 \leq {\rm Re}\,z \leq 1/2, \, {\rm Im}\, z = 1/2\},\\
S_{7}& = {\mathcal M}\cap \{-1/2 \leq {\rm Re}\,z \leq 0, \, |z|^2 + y^2 = 1\},\\
S_{8}& = {\mathcal M}\cap \{0 \leq {\rm Re}\,z \leq 1/2, \, |z|^2 + y^2 = 1\}.
\end{split}
\nonumber
\end{equation}
The side-pairing is given as follows :
\begin{equation}
h\cdot : S_{1} \to S_{2},
\label{S3Sidepairing1}
\end{equation}
\begin{equation}
g_3\cdot : S_{3} \to S_{4},
\label{S3Sidepairing2}
\end{equation}
\begin{equation}
g_3\cdot g_2\cdot g_3\cdot : S_{5} \to S_{6},
\label{S3Sidepairing3}
\end{equation}
\begin{equation}
g_1\cdot : S_{7} \to S_{8}.
\label{S3Sidepairing4}
\end{equation}

\subsubsection{Singular points} 
We next compute fixed points for $PSL(2,{\bf Z}[i])$.
We put
\begin{equation}
L_{ij} = S_i\cap S_j.
\nonumber
\end{equation}
For a set $A$ in $\mathcal M$, its isotropy group, ${\mathcal I}\,(A)$, is defined by
\begin{equation}
{\mathcal I}\,(A) = \{\gamma \in PSL(2,{\bf Z}[i])\, ; \, \gamma\cdot\zeta = \zeta, \ \forall \zeta \in A\}.
\nonumber
\end{equation}
If $\zeta = z + yj \in \mathcal M$ is a fixed point of $1 \neq \gamma \in PSL(2,{\bf Z}[i])$, we have by (\ref{S3fixedpoints})
\begin{equation}
\left\{
\begin{split}
& |cz + d|^2 + |cy|^2 = 1, \\
& (az + b)\overline{(cz + d)} + a\overline{c}y^2 = z,  \\
& ad - bc = 1.
\end{split}
\right.
\label{S3Equations}
\end{equation}
From (\ref{S3Equations}), we have $|c|^2y^2 \leq 1$. On the other hand, $y^2 \geq 1/2$ by Lemma 3.3 (2). Therefore, $|c|^2 = 0, 1, 2$. Let us compute these 3 cases separately.

\medskip
\noindent
Case 1) If $|c|^2 = 0$, we have by (\ref{S3Equations})
$$
|d|=1, \quad a\overline{d}z + b\overline{d} = z, \quad ad = 1.
$$ 
Since $\gamma \neq 1$, $d \neq \pm 1$. By the definition of $PSL(2,{\bf Z}[i])$, we have only to consider the case $d = i$. By a direct computation we have the following results :

\medskip
\begin{center}
\begin{tabular}{|c|c|c|c|c|c|}
\hline
$a$ & $b$ & $c$ & $d$ & Singular point $\zeta$ & Singular set \\
\hline
$-i$ & 0 & 0 & $i$ & $z = 0, \ y \geq 1$ & $L_{34}$ 
\\
\hline
$-i$ & $- 1$ & 0& $i$ & $z = i/2, \ y \geq {\sqrt3}/{2}$ & $L_{56}$   \\
\hline 
$-i$ & $-i$& $0$ & $i$ & $z = -1/2, \ y \geq \sqrt3/2$ & $L_{13}$ \\
\hline
$-i$ &$ i$  &$0$ & $i$& $z = 1/2, \ y \geq \sqrt3/2$& $L_{24}$ \\
\hline
$-i$ & $-1+i$ & $0$ & $i$ & $z = 1/2 + i/2, \ y \geq\sqrt2/2$& $L_{26}$ \\
\hline
$-i$ &$-1-i$ &$0$ &$i$ & $z = -1/2 + i/2, \ y \geq\sqrt2/2$ & $L_{15}$\\
\hline
\end{tabular}
\end{center}
\medskip
\begin{center}
{\bf Table $A$}
\end{center}

\medskip
\noindent
Case 2) If $|c|^2 = 1$, we need to consider 2 cases : $c = 1, i$. Let
\begin{equation}
S(c,d) = \{z + yj\, ; \, |cz + d|^2 + |cy|^2 = 1\}.
\nonumber
\end{equation}  
When $|c| =1$, $S(c,d)\cap\mathcal M \neq \emptyset$ only if $|d|^2 \leq 2$. 
 In fact, since $d \in {\bf Z}[i]$, $|d|^2 > 2$ implies $|d|^2 \geq 4$. Therefore, a simple geometry shows $1 \geq |cz + d|^2 \geq (|d|-|z|)^2 \geq (2 - 1/\sqrt{2})^2$, which is a contradiction.
If $|d|^2 = 2$, then $d = \pm (1 \pm i)$, and for $c = 1$ or $i$, we see that $\mathcal M \cap S(c,\pm (1 \pm i)) = \emptyset$ again by a geometric observation. Therefore $d$ must be $\pm 1, \pm i, 0$. Let us note that if $|\alpha|=1$, $|z+\alpha|^2 + y^2=1$ and $|z|^2 + y^2 > 1$, then ${\rm Re}\,(\overline{\alpha}z) < - 1/2$. Hence
$$
S(1,\alpha)\cap\mathcal M \subset \{|z|^2 + y^2 = 1\},
$$
if $\alpha \in {\bf Z}[i]$, $|\alpha|=1$. This implies 
\begin{equation}
S(1,\alpha)\cap\mathcal M = \{|z|^2+y^2=1,  \ {\rm Re}\,(\overline{\alpha}z) = - 1/2, \ |{\rm Re}\, z| \leq 1/2, \ 0 \leq {\rm Im}\,z \leq 1/2\}.
\label{S2S1alpha}
\end{equation}

\medskip
\noindent
(2-1) If $c = 1, d = 1$, we have by (\ref{S2S1alpha})
\begin{equation}
S(1,1)\cap\mathcal M = \{|z|^2 + y^2 = 1, \, {\rm Re}\,z = -1/2,\ 0\leq {\rm Im}\,z \leq 1/2\}.
\nonumber
\end{equation}
The last 2 equalities of (\ref{S3Equations}) and $|z|^2 + y^2 = 1$ imply 
$2b(1 + {\rm Re}\,z) = - 1$. Hence
\begin{equation}
\quad \zeta \in L_{17}, \quad a = 0, \quad b = -1.
\nonumber
\end{equation}

\medskip
\noindent
(2-2)
If $c = 1, d = -1$, we have, similarly
\begin{equation}
S(1,-1)\cap\mathcal M = \{|z|^2 + y^2 = 1, \, {\rm Re}\,z = 1/2,\ 0\leq {\rm Im}\,z \leq 1/2\}.
\nonumber
\end{equation}
Hence we have
\begin{equation}
\zeta \in L_{28}, \quad a = 0, \quad b = -1.
\nonumber
\end{equation}

\medskip
\noindent
(2-3)
If $c = 1$, $d = i$, then $S(1,i)\cap \mathcal M = \emptyset$.

\medskip
\noindent
(2-4)
If $c = 1$, $d = -i$,
$$
S(1,-i)\cap\mathcal M = \{|z|^2 + y^2 = 1, \ |{\rm Re}\, z| \leq 1/2, \, {\rm Im}\, z = 1/2\}.
$$
Let $z = x_1 + i/2, |x_1| \leq 1/2$. 
The last 2 equalities of (\ref{S3Equations}) and $|z|^2 + y^2 = 1$ imply $b = -1-2iz = - 2ix_1$. Since $b \in {\bf Z}[i]$, we have $x_1 = \pm 1/2, 0$. Then 
\begin{equation}
(z,a,b) = (-1/2+i/2,-1+i,i), \ (i/2,i,0), \ (1/2+i/2,1+i,-i).
\nonumber
\end{equation}

\medskip
\noindent
(2-5)
If $c = 1$, $d = 0$, then 
$$
S(1,0)\cap\mathcal M = \{|z|^2 + y^2 = 1, \ |{\rm Re}\,z| \leq 1/2, \ 
0 \leq {\rm Im}\,z \leq 1/2\}.
$$
Computing as above, we have
\begin{equation}
\begin{split}
\zeta \in L_{17}, \quad a = -1, \quad b = - 1, \\
\zeta \in L_{78}, \quad a=0, \quad b =-1, \\
\zeta \in L_{28}, \quad a = 1, \quad b = -1.
\end{split}
\nonumber
\end{equation}
(2-6)
If $c = i$, $d = 1$, then $S(i,1) = S(1,-i)$, Therefore
$$
S(i,1)\cap\mathcal M = \{|z|^2 + y^2 = 1, \ |{\rm Re}\, z| \leq 1/2, \, {\rm Im}\, z = 1/2\}.
$$
In this case, (\ref{S3Equations}) has no solution. 

\medskip
\noindent
(2-7)
If $c = i$, $d = -1$, then $S(i,-1) = S(1,i)$. Therefore by (2-3), $S(i,-1)\cap \mathcal M = \emptyset$.

\medskip
\noindent
(2-8)
If $c = i, d = i$, then $S(i,i) = S(1,1)$. Therefore
$$
S(i,i)\cap\mathcal M = \{|z|^2 + y^2 = 1, \, {\rm Re}\,z = -1/2,\ 0\leq {\rm Im}\,z \leq 1/2\}.
$$
Letting $z = -1/2 + ix_2$, we have, as above, $b = i + 2iz = - 2x_2$. Since $b \in {\bf Z}[i]$, $x_2 = 0, 1/2$. We thus have
$$
(z,a,b) = (-1/2, -i, 0), \ (-1/2 + i/2, - 1 - i, -1).
$$
(2-9)
If $c = i$, $d = -i$, then $S(i,-i) = S(1,-1)$. Theorefore,
$$
S(i,-i)\cap \mathcal M = \{|z|^2 + y^2 = 1, \, {\rm Re}\,z = 1/2,\ 0\leq {\rm Im}\,z \leq 1/2\},
$$
$$
(z,a,b) = (1/2, i, 0), \ (1/2 + i/2, -1+i,1).
$$
(2-10)
If $c = i$, $d = 0$, we have
$$
S(i,0)\cap\mathcal M = \{|z|^2 + y^2 = 1, \ |{\rm Re}\,z| \leq 1/2, \ 0 \leq {\rm Im}\,z \leq 1/2\},
$$
\begin{equation}
\begin{split}
\zeta \in L_{37}\cup L_{48}, \quad a = 0, \quad b = i, \\
\zeta \in L_{57}\cup L_{68}, \quad a = -1, \quad b = i.
\end{split}
\nonumber
\end{equation}

Summarizing, we have the following table.

\begin{center}
\begin{tabular}{|c|c|c|c|c|c|}
\hline
$a$ & $b$ & $c$ & $d$ & Singular point $\zeta$ & Singular set \\
\hline
$0$ & $-1$ & $1$ & $1$ & $|z|^2 + y^2 = 1, \, {\rm Re}\,z = -1/2,\ 0\leq {\rm Im}\,z \leq 1/2$ & $L_{17}$ 
\\
\hline
$0$ & $-1$ & $1$ & $-1$ & $|z|^2 + y^2 = 1, \, {\rm Re}\,z = 1/2,\ 0\leq {\rm Im}\,z \leq 1/2$ & $L_{28}$ 
\\
\hline
$-1+i$ & $i$ & $1$& $-i$ & $z = -1/2 + i/2, \ y =  1/\sqrt2$ & $L_{17}\cap L_{57}\cap L_{15}$   \\
\hline 
$i$ & $0$ & $1$& $-i$ & $z = i/2, \ y = \sqrt3/2$ & $L_{57}\cap L_{68}\cap L_{56}$ \\
\hline
$1 + i$ & $-i$ & $1$ & $-i$ & $z = 1/2 + i/2, \ y = 1/\sqrt2$ & $L_{28}\cap L_{68}\cap L_{26}$ \\
 \hline
 $-1$ & $-1$ & $1$ & $0$ & $|z|^2 + y^2 = 1, \ {\rm Re}\,z = -1/2, \ 0 \leq {\rm Im}\,z \leq 1/2$ & $L_{17}$\\
 \hline
 $0$ & $-1$ & $1$ & $0$ & $|z|^2+y^2 = 1, \ {\rm Re}\,z = 0, \ 0 \leq {\rm Im}\,z \leq 1/2$ & $L_{78}$\\
 \hline
 $1$ & $-1$ & $1$ & $0$ & $|z|^2 + y^2 = 1, \ {\rm Re}\,z = 1/2, \ 0 \leq {\rm Im}\,z \leq 1/2$ & $L_{28}$\\
 \hline
 $-i$ & $0$ & $i$ & $i$ & $z = -1/2, \ y = \sqrt3/2$ & $L_{17}\cap L_{37}\cap L_{13}$\\
 \hline 
 $-1-i$ & $-1$ & $i$ & $i$ & $z = -1/2 + i/2, \ y = 1/\sqrt2$ & $L_{17}\cap L_{57}\cap L_{15}$\\
 \hline
 $i$ & $0$ & $i$ & $-i$ & $z = 1/2, \ y = \sqrt3/2$ & $L_{28}\cap L_{48}\cap L_{24}$ \\
 \hline
 $-1+i$ & $1$ & $i$ & $-i$ & $z = 1/2 + i/2, \ y = 1/\sqrt{2}$ & $L_{28}\cap L_{68}\cap L_{26}$ \\
\hline
$0$ & $i$ & $i$ & $0$ & $ |z|^2 + y^2 = 1, \ |{\rm Re}\,z| \leq 1/2, \ {\rm Im}\,z = 0$ & $L_{37}\cup L_{48}$\\
 \hline
 $-1$ & $i$ & $i$ & $0$ & $|z|^2+y^2=1, \ |{\rm Re}\,z| \leq 1/2, \ 
 {\rm Im}\,z = 1/2$ & $L_{57}\cup L_{68}$\\
 \hline
\end{tabular}
\end{center}
\begin{center}
{\bf Table B}
\end{center}

\medskip
\noindent
Case 3) Since $y^2 \geq 1/2$, when $|c|^2 = 2$, we have $y^2 = 1/2$, $cd + z=0$, which inturn implies $d\neq 0$, as $ad-bc=1$. Therefore, we have $|z|^2 = 1/2$, which implies
$$
\zeta = \pm \frac{1}{2} + \frac{1}{2}i + \frac{1}{\sqrt2}j, \quad 
c = 1 \pm i.
$$
Computing $a, b, d$ by the formulas
\begin{equation}
cz + d = 0, \quad a\overline{c} = 2z, \quad
 ad - bc = 1,
\nonumber
\end{equation}
we get the following table.

\bigskip
\begin{center}
\begin{tabular}{|c|c|c|c|c|c|}
\hline
$a$ & $b$ & $c$ & $d$ & $\zeta$ &Singular set \\
\hline
$i$ & $0$ & $1 + i$ & $-i$ & $z = 1/2 + i/2, \ y = 1/\sqrt2$ & $L_{28}\cap L_{68}\cap L_{26}$\\
\hline
$-1$ & $-1+i$ & $1 + i$ & $1$ & $z = -1/2 + i/2, \ y=1/\sqrt2$& $L_{17}\cap L_{57}\cap L_{15}$\\
\hline
$1$ & $-1-i$&$1 -i$ &$-1$ &$z = 1/2 + i/2, \ y = 1/\sqrt2$ &$L_{28}\cap L_{68}\cap L_{26}$\\
\hline 
$i$ &$0$ & $1-i$& $-i$ &$z = -1/2 + i/2, \ y = 1/\sqrt2$& $L_{17}\cap L_{57}\cap L_{15}$\\
 \hline
\end{tabular}
\end{center}
\begin{center}
{\bf Table C}
\end{center}
\bigskip
By virtue of Tables A, B, C, the set of singular points is split into 13 parts :
\begin{equation}
\begin{split}
L_{13} &= \{z = -1/2, \ y \geq \sqrt3/2\}, \\
 L_{34} & = \{z = 0, \ y \geq 1\},\\
L_{24} &= \{z = 1/2, \quad y \geq \sqrt3/2\}, \\
L_{15} &= \{z = -1/2 + i/2, \quad y \geq 1/\sqrt2\}, \\
 L_{56} & = \{z = i/2, \ y\geq \sqrt3/2\}, \\
L_{26} &= \{z = 1/2 + i/2, \quad y \geq 1/\sqrt2\}, \\
 L_{17} &= \{|z|^2+y^2 = 1, \ {\rm Re}\,z = -1/2, \ 0 \leq {\rm Im}\,z \leq 1/2\}, \\
 L_{78} & = \{|z|^2+y^2 = 1, \ {\rm Re}\,z = 0, \ 0 \leq {\rm Im}\,z \leq 1/2\}, \\
 L_{28} &= \{|z|^2+y^2 = 1, \ {\rm Re}\,z = 1/2, \ 0 \leq {\rm Im}\,z \leq 1/2\},\\
 L_{37} &= \{|z|^2+y^2 = 1, \ -1/2 \leq {\rm Re}\,z \leq 0, \ {\rm Im}\,z = 0\},\\
 L_{48}& = \{|z|^2+y^2 = 1, \ 0 \leq {\rm Re}\,z \leq 1/2, \ {\rm Im}\,z = 0\},\\
L_{57} &= \{|z|^2+y^2 = 1, \ -1/2 \leq {\rm Re}\,z \leq 0, \ {\rm Im}\,z = 1/2\},\\
 L_{68} &= \{|z|^2+y^2 = 1, \ 0 \leq {\rm Re}\,z \leq 1/2, \ {\rm Im}\,z = 1/2\}.
\end{split}
\nonumber
\end{equation}
By the side-pairing (\ref{S3Sidepairing1}) $\sim$ (\ref{S3Sidepairing4}), we have the following identification :
\begin{equation}
L_{13} = L_{24}, \quad L_{15} = L_{26}, \quad
L_{17} = L_{28}, \quad L_{37} = L_{48}, \quad L_{57} = L_{68}.
\label{S2Lidentify}
\end{equation}
With this in mind, we put
\begin{equation}
\begin{split}
\mathcal L_1 &= L_{13}, \quad \mathcal L_2 = L_{15}, \quad
\mathcal L_3 = L_{34}, \quad \mathcal L_4 = L_{56}, \\
\mathcal L_5& = L_{17}, \quad \mathcal L_6 = L_{78}, \quad
\mathcal L_7 = L_{37}, \quad 
\mathcal L_{8} = L_{57}.
\end{split}
\nonumber
\end{equation}
The isotropy groups for $\mathcal L_n$ are as follows :
\begin{equation}
{\mathcal I}\,(\mathcal L_1) = \left\langle
\left(
\begin{array}{cc}
-i & -i\\
0 & i
\end{array}
\right)
\right\rangle = {\bf Z}_2,
\nonumber
\end{equation}
\begin{equation}
{\mathcal I}\,(\mathcal L_2) = \left\langle
\left(
\begin{array}{cc}
-i & -1-i\\
0 & i
\end{array}
\right)
\right\rangle = {\bf Z}_2,
\nonumber
\end{equation}
\begin{equation}
{\mathcal I}\,(\mathcal L_3) = \left\langle
\left(
\begin{array}{cc}
-i & 0\\
0 & i
\end{array}
\right)
\right\rangle = {\bf Z}_2,
\nonumber
\end{equation}
\begin{equation}
{\mathcal I}\,(\mathcal L_4) = \left\langle
\left(
\begin{array}{cc}
-i & - 1\\
0 & i
\end{array}
\right)
\right\rangle = {\bf Z}_2,
\nonumber
\end{equation}
\begin{equation}
{\mathcal I}\,(\mathcal L_5) = \left\langle
\left(
\begin{array}{cc}
0 & -1\\
1 & 1
\end{array}
\right)
\right\rangle = {\bf Z}_3,
\nonumber
\end{equation}
\begin{equation}
{\mathcal I}\,(\mathcal L_6) = \left\langle
\left(
\begin{array}{cc}
0 & -1\\
1 & 0
\end{array}
\right)
\right\rangle = {\bf Z}_2,
\nonumber
\end{equation}
\begin{equation}
{\mathcal I}\,(\mathcal L_7) = \left\langle
\left(
\begin{array}{cc}
0 & i\\
i & 0
\end{array}
\right)
\right\rangle = {\bf Z}_2,
\nonumber
\end{equation}
\begin{equation}
{\mathcal I}\,(\mathcal L_8) = \left\langle
\left(
\begin{array}{cc}
-1 & i\\
i & 0
\end{array}
\right)
\right\rangle = {\bf Z}_3.
\nonumber
\end{equation}


\subsection{Uniformizing cover}
For $p \in \mathcal M$, we define its uniformizing cover $B(\widetilde p\,;\,\widetilde g_p)$ in 
${\bf H}^3$, where $\widetilde p = p$ as a point in ${\bf H}^3$, and $\widetilde g_p$ is the 
hyperbolic metric on ${\bf H}^3$. We give the radius $r_p$ of $B(\widetilde p\,;\,\widetilde g_p)$ 
and $G_p$.  If $p \not\in \cup_n\mathcal L_n$, we take $0 < r_p < {\rm dist}(p,\cup_n\mathcal L_n)$, 
$G_{p} = \{1\}$.
If $p$ is in the interior of $\mathcal L_n$, we take 
$0 < r_p < {\rm dist}(p,\cup_{m\neq n}\mathcal L_m)$, 
$G_p = {\rm Iso}\,(\mathcal L_n)$.
Taking account of the side-pairing, we put
\begin{equation}
P_1  = - \frac{1}{2} + \frac{\sqrt3}{2}j, \quad P_2 = - \frac{1}{2} + \frac{1}{2}i + \frac{\sqrt2}{2}j, \quad P_3 = j,  \quad
P_4 = \frac{1}{2}i + \frac{\sqrt3}{2}j.
\nonumber
\end{equation}
If $p = P_m$, we take $r_p$ small enough so that $B(\widetilde p\,;\,\widetilde g_p)$ does 
not intersect with $\mathcal L_n$ which does not meet $p$, and $G_p$ is the group generated by 
${\rm Iso}\,(\mathcal L_n)$, for all $n$ such that $\mathcal L_n$ meets $p$.
If $U_p\cap U_q \neq \emptyset$ and $\pi_p(\zeta) = \pi_q(\zeta')$, we let 
$V_{\zeta}$, $V_{\zeta'}$ be small balls centered at $\zeta$, $\zeta'$, and define 
$\psi$ as an isometry from $V_{\zeta}$ to $V_{\zeta'}$ such that $\psi(\zeta) = \zeta'$, 
the construction of which will be evident. This defines the orbifold structure on 
$\mathcal M$.


\subsection{Orbifold structure of the horizontal slice}
Since the isotropy groups for $\mathcal L_1$, $\mathcal L_2$, $\mathcal L_3$, $\mathcal L_4$, 
are the ones of the rotation around them, for any  $t > 1$, the horizontal slice 
$M_t = \mathcal M\cap\{y = t\}$ is a compact 2-dimensional orbifold with 
singular points 
$
-1/2 + tj,\  -1/2 + i/2 + tj, \ tj, \ i/2 + tj.
$
Here the covering is a disc in ${\bf R}^2$, and we use the scaled Euclidean metric instead of the 
hyperbolic metric. 
Note that the orbifold structure of $M_t$ is independent of $t$.


\section{Spectral properties of the model space} 

In this section, we study spectral properties of the Laplace-Beltrami operator $H_{free}$ for the 
space which models an end of our orbifold $\mathcal M$. Namely, we take 
$\mathcal M_{free}= M_{free} \times (0, \infty)$ with metric
\begin{equation} \label{3.09.1}
ds^2= \frac{(dy)^2+h_{free}(x, dx)}{y^2},
\end{equation}
where $M_{free}$ is a compact $(n-1)-$dimensional Riemannian orbifold with metric $h_{free}$,
cf. (\ref{S1ds2regular}). If $U_j,\, j=1, \dots, m,$ is a finite uniformizing covering 
of $M_{free}$ by coordinate
patches with $U_j= \widetilde U_j/ G_j$, see \S 1.4, we use
the partition of unity $\chi_j(x, y)=\phi_j(x),\, j=1, \dots,m$, where $\phi_j$ 
form a smooth partition of unity on $M_{free}$, cf. \S 1.4.


\subsection{Besov type spaces}
We first discuss a Besov type space suitable to treat the resolvent of $H_0$. It is a hyperbolic
analog of the Besov type space
introduced by Agmon-H{\"o}rmander in the case of the Euclidean space.

Let $\bf H$ be a Hilbert space endowed with norm $\|\cdot\|_{\bf H}$.
We decompose $(0,\infty)$ into $(0,\infty) = \cup_{k\in {\bf Z}}I_k$, where
\begin{equation}
I_k = \left\{
\begin{array}{cc}
\big(\exp(e^{k-1}),\exp(e^k)\big], & k \geq 1, \\
\big(e^{-1},e\big], & k = 0, \\
\big(\exp(-e^{|k|}),\exp(-e^{|k|-1})\big], & k \leq -1. 
\end{array}
\right.
\nonumber
\end{equation}
Let $\mathcal B$ be the Banach space of ${\bf H}$-valued function on $(0,\infty)$ such that
\begin{equation}
\|f\|_{\mathcal B} = \sum_{k\in{\bf Z}}e^{|k|/2}
\left(\int_{I_k}\|f(y)\|_{\bf H}^2\, \frac{dy}{y^n}\right)^{1/2} < \infty.
\label{S3Bspace}
\end{equation}
The dual space of $\mathcal B$ is identified with the space equipped with norm 
\begin{equation}
\|u\|_{\mathcal B^{\ast}} = \left(\sup_{R>e}\frac{1}{\log R}\int_{\frac{1}{R}<y<R}\|u(y)\|_{\bf H}^2\, \frac{dy}{y^n}\right)^{1/2} < \infty.
\label{S3Bast}
\end{equation}
Therefore, we have
\begin{equation}
\left|\int_0^{\infty}(f(y),u(y))_{\bf H}\frac{dy}{y^n}\right| \leq C\|f\|_{\mathcal B}\|u\|_{\mathcal B^{\ast}},
\label{S5L2andBesov}
\end{equation}
with a constant $C > 0$.
We also use the following weighted $L^2$-space. For $s \in {\bf R}$, 
\begin{equation}
L^{2,s} \ni u \Longleftrightarrow \|u\|_{s}^2 = \int_0^{\infty}
(1 + |\log y|)^{2s}\|u(y)\|_{\bf H}^2\,\frac{dy}{y^n} < \infty.
\label{S3L2s}
\end{equation}
For $s > 1/2$, the following inclusion relation holds :
\begin{equation}
L^{2,s} \subset \mathcal B \subset L^{2,1/2} \subset L^2 \subset L^{2,-1/2} \subset \mathcal B^{\ast} \subset L^{2,-s}.
\label{S3Inclusion}
\end{equation}


\begin{lemma} \label{L:5.1} The following two assertions are equivalent.
\begin{equation}
\lim_{R\to\infty}\frac{1}{\log R}\int_{\frac{1}{R}<y<R}\|u(y)\|_{\bf H}^2\,\frac{dy}{y^n} = 0.
\end{equation}
\begin{equation}
\lim_{R\to\infty}\frac{1}{\log R}\int\rho\Big(\frac{\log y}{\log R}\Big)\|u(y)\|_{\bf H}^2\,\frac{dy}{y^n} = 0, \quad \forall \rho \in C_0^{\infty}((0,\infty)).
\end{equation}
\end{lemma}

The proof of the above results are given in \cite{IsKu09}, Chap. 1, \S 2.


\subsection{A-priori estimates}
Let $(M_{free}, h_{free})$ be a compact $(n-1)-$dimensional Riemannian orbifold with Laplace 
operator $\Delta_{free}$, and $dV_{free}$ the volume element of $M_{free}$. The inner product of $L^2(M_{free})$ is denoted by $(f,g)_{L^2(M_{free})}$. Let
$\mathcal M_{free} = M_{free}\times(0,\infty)$, which is an $n$-dimensional Riemannian orbifold with metric $y^{-2}((dy)^2 + h_{free}(x,dx))$. We  put 
\begin{equation}
H_{free} = - y^2(\partial_y^2 + \Delta_{free}) + (n-2)y\partial_y - \frac{(n-1)^2}{4},
\label{S5H0}
\end{equation}
which, restricted on $C_0^{\infty}(\mathcal M_{free})$, is essentially self-adjoint in
\begin{equation}
{\mathcal H}_{free} = L^2\Big({\mathcal M_{free}};\frac{dV_hdy}{y^n}\Big) 
 = L^2\Big((0,\infty);L^2(M_{free});\frac{dy}{y^n}\Big).
\label{S5spaceH0}
\end{equation}
Let $0 = \lambda_0 \leq \lambda_1 \leq \lambda_2 \leq \cdots$ be the eigenvalues of $- \Delta_{free}$ 
with complete orthonormal system of eigenvectors 
$|M_{free}|^{-1/2} = \varphi_0, \varphi_1, \varphi_2, \cdots$, $|M_{free}|$ being the volume 
of $M_{free}$. 
We use the following notation:
\begin{equation}
D_y = y\partial_y, \quad D_x = y\sqrt{-\Delta_{free}}.
\label{S5DyDx}
\end{equation}
Let $\mathcal B$ and $\mathcal B^{\ast}$ be as in Subsection 4.1 with ${\bf H}$ replaced 
by $L^2(M_{free})$.


\begin{lemma}\label{Aprioriestimate}
(1) If $u \in \mathcal B^{\ast}$ satisfies $(H_{free} - z)u = f \in \mathcal B^{\ast}$ with $z \in {\bf C}$, we have
$$
\|D_yu\|_{\mathcal B^{\ast}} + \|D_xu\|_{\mathcal B^{\ast}} \leq 
C(\|u\|_{\mathcal B^{\ast}} + \|f\|_{\mathcal B^{\ast}}).
$$
(2) Furthermore, if
$$
\lim_{R\to\infty}\frac{1}{\log R}\int_{1/R}^R\left[\|u(y)\|^2_{L^2(M_{free})} + 
\|f(y)\|^2_{L^2(M_{free})}\right]\frac{dy}{y^n} = 0
$$
holds, we have
$$
\lim_{R\to\infty}\frac{1}{\log R}\int_{1/R}^R\left[\|D_yu(y)\|^2_{L^2(M_{free})} + 
\|D_xu(y)\|^2_{L^2(M_{free})}\right]\frac{dy}{y^n} = 0.
$$
where $\|\cdot\|$ denotes the norm of $L^2(\mathcal M_{free})$.\\
(3) The assertion (2) also holds with ${\rm lim}$ replaced by ${\rm lim \; inf}$. \\
(4) If $(H_{free} -z)u = f$ holds and $u, f \in L^{2,s}$ for some $s \geq 0$. Then 
$$
\|D_yu\|_s + \|D_xu\|_s \leq C_s(\|u\|_s + \|f\|_s),
$$

In the above estimates in (1) and (4), the constant $C$ is independent of $z$ when $z$ varies over a bounded set in ${\bf C}$. \\
\noindent
(5) If $f \in L^{2,s}$ for some $s \geq 0$, $u = R_{free}(z)f$ $(z \not\in {\bf R})$ satisfies
$$
\|u\|_s + \|D_yu\|_s + \|D_xu\|_s \leq C_{s,z}\|f\|_s,
$$
where the constant $C_{s,z}$ is independent of $z$ when $z$ varies over a compact set in ${\bf C} \setminus{\bf R}$.
\end{lemma}

Proof. 
Let 
$$
u_m(y) = (u(y),\varphi_m)_{L^2(M_{free})}, \quad
f_m(y) = (f(y),\varphi_m)_{L^2(M_{free})}.
$$
Then $u_m$ satisfies
\begin{equation}
(- D_y^2 + (n-1)D_y + y^2\lambda_m - z)u_m = f_m.
\label{S5Equationuj}
\end{equation}
We take $\chi(t) \in C_0^{\infty}({\bf R})$ such that $\chi(t) = 1$ $(|t| < 1)$, $\chi(t) = 0$ $(|t| > 2)$, and put
$$
\chi_R(y) = \chi\Big(\frac{\log y}{\log R}\Big).
$$
By integration by parts, we have
$$
\|\chi_RD_yu_m\|^2 + \lambda_m\|y\chi_Ru_m\|^2 = 
z\|\chi_Ru_m\|^2 - (D_yu_m,(D_y\chi_R^2)u_m) + (f_m,\chi_R^2u_m),
$$
which implies
\begin{equation}
\begin{split}
& \|\chi_RD_yu_m\|^2 + \lambda_m\|y\chi_Ru_m\|^2 \\
& \leq |z|\|\chi_Ru_m\|^2 + \frac{2}{\log R}|(\chi_RD_yu_m,\chi'\big(\frac{\log y}{\log R}\big)u_m)| + |(f_m,\chi_R^2u_m)|.
\end{split}
\nonumber
\end{equation}
We then have for $R > e^2$
\begin{equation}
\|\chi_RD_yu_m\|^2 + \lambda_m\|y\chi_Ru_m\|^2 \leq
C_z\Big(\|\chi_Ru_m\|^2 + \|\chi'\big(\frac{\log y}{\log R}\big)u_m\|^2 + \|\chi_Rf_m\|^2\Big).
\label{S5chiRDyu1}
\end{equation}
We sum up these inequalities with respect to $m$ and divide it by $\log R$. 
Taking the supremum with respect to $R$, 
we get the assertion (1).

Letting $R \to \infty$  and using Lemma \ref{L:5.1}, we get (2) and (3). 

Letting $R \to \infty$ in (\ref{S5chiRDyu1}) and then summing up the resulting inequalities
with respect to $m$, 
we obtain (4) for $s = 0$. Let $0 \leq s \leq 1$, and put $v = (1 + |\log y|^2)^{s/2}u$. Then $v$ satisfies 
$$
(H_{free} - z)v = (1 + |\log y|^2)^{s/2}f + Pu,
$$
where $P$ is a 1st order differential operator with respect to $D_y, D_x$, 
whose coefficients decay like $O(1 + |\log y|)^{s-1}$, so that, by the proven result with $s=0$,
$P u \in L^2$ when $s \leq 1$.  
Applying again the result for $s = 0$, we obtain the $L^2-$estimate for $v$. Thus,
we get (4) for when $0 \leq s \leq 1$. Similarly, one can prove (4) for all $s \geq 0$. 

By the standard $L^2-$resolvent estimates, (5) follows easily from (4).
\qed


\subsection{Identity} 
In the above proof, we have seen that the 1-dimensional operator
\begin{equation}
L_{free}(\zeta) = y^2(- \partial_y^2 + \zeta^2) + (n-2)y\partial_y - \frac{(n-1)^2}{4}
\label{S5L0zeta}
\end{equation}
plays an important role. It has the following Green kernel
\begin{equation}
G_{free}(y,y';\zeta,\nu) = 
\left\{
\begin{split}
(yy')^{(n-1)/2}K_{\nu}(\zeta y)I_{\nu}(\zeta y'), \quad y > y' > 0, \\
(yy')^{(n-1)/2}I_{\nu}(\zeta y)K_{\nu}(\zeta y'), \quad y' > y > 0,
\end{split}
\right.
\label{S5G0yy'}
\end{equation}
where $I_{\nu}(z)$, $K_{\nu}(z)$ are the modified Bessel functions.
Namely,
$$
(L_{free}(\zeta) + \nu^2)^{-1}f = \int_0^{\infty}G_{free}(y,y';\zeta,\nu)f(y')\frac{dy'}{(y')^n}, 
\quad \forall f \in C_0^{\infty}((0,\infty)).
$$
Using the asymptotic expansion of the modified Bessel functions, we observe the behavior 
of $v_{\pm} = (L_{free}(\zeta) - k^2 \mp i0)^{-1}f,$ $k > 0,$ to see that
$$
v_{\pm} \sim C_{\pm}(k)y^{(n-1)/2 \mp ik}, \quad y \to 0.
$$
 Therefore, we infer
$$
\left(y\partial_y - \big(\frac{n-1}{2} \mp ik\big)\right)v_{\pm} = o(y^{(n-1)/2}), \quad y \to 0.
$$
This suggests the importance of the term $\left(y\partial_y - \big(\frac{n-1}{2} \mp i\sqrt{z}\big)\right)(H_{free} - z)^{-1}f$ to derive the estimates of the resolvent. We put
\begin{equation}
\sigma_{\pm} = \sigma_{\pm}(z) = \frac{n-1}{2} \mp i\sqrt{z},
\label{S5sigmapm}
\end{equation}
where for $z = re^{i\theta}$, $r > 0$, $- \pi < \theta < \pi$, we take the branch of $\sqrt{z}$ as $\sqrt{r}e^{i\theta/2}$. 

In the following Lemmas of this section, $(\;,\,)_{0}$ and $\|\cdot\|_{0}$ denote the inner product and the norm of $L^2(M_{free})$. Our procedure for the resolvent estimates leans over two identities, those in Lemma 4.3 and (\ref{eq:Imwyyn-1inte}).


\begin{lemma} Suppose $u$ satisfies $(H_{free} - z)u = f$, and let $w_{\pm} = (D_y - \sigma_{\pm})u$.  Let $\varphi(y) \in C^1((0,\infty);{\bf R})$ and $0 <a < b < \infty$. Then we have
\begin{eqnarray*}
& &  \int_a^b(D_y\varphi + 2\varphi)\|D_xu\|_0^2\frac{dy}{y^n} + 
\left[\frac{\varphi(\|w_{\pm}\|_0^2 - \|D_xu\|_0^2)}{y^{n-1}}\right]_{y=a}^{y=b}\\
& =& \mp 2\,{\rm Im}\,\sqrt{z}\int_a^b\varphi\left(\|w_{\pm}\|_0^2 + 
\|D_xu\|_0^2\right)\frac{dy}{y^n} \\
& & + \int_a^b(D_y\varphi)\|w_{\pm}\|_0^2\frac{dy}{y^n} 
- 2\,{\rm Re}\int_a^b\varphi(f,w_{\pm})_0\frac{dy}{y^n}.
\end{eqnarray*}
\end{lemma}

Proof. 
We rewrite the equation $(H_{free} - z)u = f$ as 
\begin{equation}
D_y(D_y - \sigma_{\pm})u = \sigma_{\mp}(D_y - \sigma_{\pm})u - y^2\Delta_{free}u - f.
\label{S5rewrite1}
\end{equation}  
Taking the inner product of (\ref{S5rewrite1}) and $\varphi w_{\pm}$, we have
\begin{equation}
\begin{split}
& \int_a^b\varphi(D_yw_{\pm},w_{\pm})_0\frac{dy}{y^n} \\
= 
&\  \sigma_{\mp}\int_a^b\varphi\|w_{\pm}\|_0^2\frac{dy}{y^n} 
-\int_a^b\varphi(y^2\Delta_{free}u,w_{\pm})_0\frac{dy}{y^n} - 
\int_a^b\varphi(f,w_{\pm})_0\frac{dy}{y^n}.
\end{split}
\label{C2S2ProofLemma2.5}
\end{equation}
Take the real part.
By integration by parts, the left-hand side is equal to
\begin{equation}
\begin{split}
& {\rm Re}\int_a^b\varphi(D_yw_{\pm},w_{\pm})_0\frac{dy}{y^n} \\ 
= & \ 
\left[\frac{\varphi\|w_{\pm}\|_0^2}{2y^{n-1}}\right]_{y=a}^{y=b} 
- \frac{1}{2}\int_a^b(D_y\varphi)\|w_{\pm}\|_0^2\frac{dy}{y^n} + 
\frac{n-1}{2}\int_a^b\varphi\|w_{\pm}\|_0^2\frac{dy}{y^n}.
\end{split}
\label{S5identity1}
\end{equation}
Here, let  us note that using
$$                                    
(-y^2\Delta_{free}u,D_yu)_0 = (v,D_yv)_0 - \|v\|_0^2, \quad v = y\sqrt{-\Delta_{free}}u,
$$
we have
\begin{eqnarray*}
& &- {\rm Re}\int_a^b\varphi\big(y^2\Delta_{free}u,w_{\pm}\big)_0\frac{dy}{y^n} 
\\
& =& \left[\frac{\varphi\|D_xu\|_0^2}{2y^{n-1}}\right]_{y=a}^{y=b} 
- \frac{1}{2}\int_a^b(D_y\varphi)\|D_xu\|_0^2\frac{dy}{y^n} 
+ \left(\dfrac{n-3}{2} - {\rm Re}\,\sigma_{\pm}\right)\int_a^b\varphi\|D_xu\|_0^2\frac{dy}{y^n}.
\end{eqnarray*}
Apply this to the 2nd term of the right-hand side of (\ref{C2S2ProofLemma2.5}).  We then have
\begin{equation}
\begin{array}{rcl}
& & {\displaystyle {\rm Re}\int_a^b\varphi(D_yw_{\pm},w_{\pm})_0\frac{dy}{y^n}} \\
& =& {\displaystyle ({\rm Re}\, \sigma_{\mp})\int_a^b\varphi\|w_{\pm}\|_0^2\frac{dy}{y^n} 
-{\rm Re}\int_a^b\varphi(y^2\Delta_{free}u,w_{\pm})_0\frac{dy}{y^n} 
- {\rm Re}\int_a^b\varphi(f,w_{\pm})_0\frac{dy}{y^n}} \\
&= &{\displaystyle 
\left(\frac{n-1}{2} \mp {\rm Im}\,\sqrt{z}\right)\int_a^b\varphi\|w_{\pm}\|_0^2\frac{dy}{y^n}
+ \left[\frac{\varphi\|D_xu\|_0^2}{2y^{n-1}}\right]_{y=a}^{y=b} } \\
& &{\displaystyle  - \frac{1}{2}\int_a^b(D_y\varphi)\|D_xu\|_0^2\frac{dy}{y^n}
-  (1 \pm {\rm Im}\,\sqrt{z})\int_a^b\varphi
\|D_xu\|_0^2\frac{dy}{y^n} - {\rm Re}\int_a^b\varphi(f,w_{\pm})_0\frac{dy}{y^n}.}
\end{array}
\label{S5identity2}
\end{equation}
Equating (\ref{S5identity1}) and (\ref{S5identity2}), we obtain the lemma. \qed


\subsection{Uniform estimates}
We shall derive estimates of the resolvent $R_{free}(z) = (H_{free} - z)^{-1}$, when $z \in {\bf C}\setminus{\bf R}$ approaches to the real axis.


\begin{lemma} 
Let $u = R_{free}(z)f$. Let $w_{\pm} = (D_y - \sigma_{\pm})u$, and put for  $C^1 \ni \varphi \geq 0$ and constants $0 < a < b$, 
\begin{equation}
L_{\pm} =  \int_a^b\big(D_y\varphi + 2\varphi\big)\|D_x u\|_0^2\frac{dy}{y^n}    
+ \left[\frac{\varphi(\|w_{\pm}\|_0^2 - \|D_xu\|_0^2)}{y^{n-1}}\right]_{y=a}^{y=b},
\label{S5Apm}
\end{equation}
\begin{equation}
R_{\pm} = \int_{a}^b(D_y\varphi)\|w_{\pm}\|_0^2\frac{dy}{y^n}
- 2{\rm Re}\int_a^b\varphi(f,w_{\pm})_0\frac{dy}{y^n}.
\label{S5Bpm}
\end{equation}
Then we have the following inequality. 
\begin{equation}
L_+ \leq R_+, \quad L_- \geq R_-, \quad {\rm if} \quad {\rm Im}\,\sqrt z \geq 0,
\label{S5ABinequality1}
\end{equation}
\begin{equation}
L_+ \geq R_+, \quad L_- \leq R_-, \quad {\rm if} \quad {\rm Im}\,\sqrt z \leq 0,
\label{S5ABinequality2}
\end{equation}
 \end{lemma}

Proof.
Using Lemma 4.3,  $\varphi \geq 0$, and the sign of ${\rm Im}\,\sqrt{z}$, we obtain the lemma.
\qed

\bigskip
In the following, $z$ varies over the region
\begin{equation}
J_{\pm} = \{z \in {\bf C} \, ; \, a \leq {\rm Re}\,z \leq b, \ 
0 < \pm {\rm Im}\,z < 1\},
\label{S5Jplusminus}
\end{equation}
where $0 < a < b$ are arbitrarily chosen constants.


\begin{lemma}\label{LemmaDxuepsilonu+fepsion}
Let $u = R_{free}(z)f$ with $f \in \mathcal B$. Then for any $\epsilon > 0$, there exists a constant $C_{\epsilon} > 0$ such that
$$
\int_0^{\infty}\|D_xu\|_0^2\frac{dy}{y^n} \leq \epsilon\|u\|_{\mathcal B^{\ast}}^2 + 
C_{\epsilon}\|f\|_{\mathcal B}^2, \quad \forall z \in J_{\pm}.
$$
\end{lemma}
Proof. Assume that $z \in J_+$. Letting $\varphi = 1$ and using (\ref{S5ABinequality1}), we have
$$
\int_a^b\|D_xu\|_0^2\frac{dy}{y^n} \leq \left[\frac{\|D_xu\|_0^2 - \|w_{+}\|_0^2}{2y^{n-1}}\right]_{y=a}^{y=b} +\left|\int_a^b(f,w_{+})_0\frac{dy}{y^n}\right|.
$$
By Lemma \ref{Aprioriestimate} (4),  $w_{+}, D_xu \in L^2$ for $z \not\in {\bf R}$. Hence
\begin{equation}
\liminf_{y\to0}\frac{\|w_{+}\|_0^2 + \|D_xu\|_0^2}{y^{n-1}} = 0, \quad
\liminf_{y\to\infty}\frac{\|w_{+}\|_0^2 + \|D_xu\|_0^2}{y^{n-1}} = 0.
\label{S5liminf}
\end{equation}
Therefore, letting $a \to 0$ and $b \to \infty$ along suitable sequences, we have
$$
\int_0^{\infty}\|D_xu\|_0^2\frac{dy}{y^n} \leq \left|\int_0^{\infty}(f,w_{+})_0\frac{dy}{y^n}\right| \leq \epsilon\|w_{+}\|^2_{\mathcal B^{\ast}} + C_{\epsilon}\|f\|_{\mathcal B}^2.
$$
Lemma \ref{Aprioriestimate} (1) yields $\|w_{+}\|_{\mathcal B^{\ast}} \leq C(\|u\|_{\mathcal B^{\ast}} + \|f\|_{\mathcal B^{\ast}})$, which proves the lemma when $z \in J_{+}$. The case for $z \in J_-$ is proved similarly by using $w_-$. \qed


\begin{lemma}  Let $u$, $f$ be as in the previous lemma,  and $w_{\pm} = (D_y - \sigma_{\pm})u$. Then for any $\epsilon > 0$, there exists a constant $C_{\epsilon} > 0$ such that for any $y > 0$ 
\begin{equation}
\frac{\|w_{+}\|_0^2 - \|D_xu\|_0^2}{y^{n-1}} \leq \epsilon\|u\|_{\mathcal B^{\ast}}^2 +
C_{\epsilon}\|f\|_{\mathcal B}^2, \quad \forall z \in J_+,
\nonumber
\end{equation}
\begin{equation}
\frac{\|w_{-}\|_0^2 - \|D_xu\|_0^2}{y^{n-1}} \leq \epsilon\|u\|_{\mathcal B^{\ast}}^2 +
C_{\epsilon}\|f\|_{\mathcal B}^2, \quad \forall z \in J_-.
\nonumber
\end{equation}
\end{lemma}
Proof. As in the previous lemma, assume that $z \in J_+$. Letting $\varphi = 1$ and using (\ref{S5ABinequality1}), we have
$$
\frac{\|w_{+}\|_0^2 - \|D_xu\|_0^2}{y^{n-1}}\Big|_{y = b} \leq 
\frac{\|w_{+}\|_0^2  - \|D_xu\|_0^2}{y^{n-1}}\Big|_{y = a} + C\|f\|_{\mathcal B}\|w_+\|_{{\mathcal B}^{\ast}}.
$$
Letting $a \to 0$ along a suitable sequence, using (\ref{S5liminf}) and Lemma \ref{Aprioriestimate}, we obtain the lemma. \qed


\begin{lemma} Let $u$, $f$, $w_{\pm}$ be as in the previous lemma.
Then for any $\epsilon > 0$, there exists a constant $C_{\epsilon} > 0$ such that
\begin{equation}
\|w_{+}\|_{{\mathcal B}^{\ast}} \leq \epsilon\|u\|_{\mathcal B^{\ast}} + C_{\epsilon}\|f\|_{\mathcal B}, \quad \forall z \in J_{+},
\nonumber
\end{equation}
\begin{equation}
\|w_{-}\|_{{\mathcal B}^{\ast}} \leq \epsilon\|u\|_{\mathcal B^{\ast}} + C_{\epsilon}\|f\|_{\mathcal B}, \quad \forall z \in J_{-}.
\nonumber
\end{equation}
\end{lemma}

Proof. We divide the inequality in Lemma 4.6 by $y$ and integrate on $(1/R,R)$. We then use Lemma 4.5 to estimate the integral of $\|D_xu\|_0^2$, and obtain 
the lemma.
\qed


\begin{theorem} There exists a constant $C > 0$ such that
$$
\|R_{free}(z)f\|_{{\mathcal B}^{\ast}} \leq C\|f\|_{\mathcal B}, \quad 
\forall z \in J_{\pm}.
$$
\end{theorem}
Proof. We consider the case that $z \in J_+$, and put $\sqrt{z} = k + i\epsilon$ for $z \in J_{+}$. Then $\epsilon > 0$ and $k > C$ for some constant $C > 0$. Letting $w_{+} = (D_y - \sigma_{+})u$, we then have 
\begin{equation}
 {\rm Im}\,\left(D_y(w_{+},u)_0\right) = 
 {\rm Im}\,\left((n - 1 + 2ik)(w_{+},u)_0\right) - {\rm Im}\,(f,u)_0.
\label{eq:ImDywu}
\end{equation}
This is a consequence of the formula
$$
D_y(w_{+},u)_0 = (D_yw_{+},u)_0 + \|w_{+}\|_0^2 + \left(\frac{n-1}{2} + \epsilon + 
ik\right)(w_{+},u)_0
$$
and (\ref{S5rewrite1}). We integrate (\ref{eq:ImDywu}). Since
$$
 \int_a^bD_y(w_{+},u)_0\frac{dy}{y^n}  = \left[\frac{(w_{+},u)_0}{y^{n-1}}\right]_a^b + (n-1)\int_a^b(w_{+},u)_0\frac{dy}{y^n},
$$
we then have
\begin{equation}
 {\rm Im}\,\left[\frac{(w_{+},u)_0}{y^{n-1}}\right]_a^b = 2k\,{\rm Re}
 \int_a^b(w_{+},u)_0\frac{dy}{y^n} - {\rm Im}\int_a^b(f,u)_0\frac{dy}{y^n}.
 \label{eq:Imwyyn-1inte}
\end{equation}
Using $w_{+} = D_yu - \sigma_{+}u$ and integrating by parts, we have
$$
 {\rm Re}\int_a^b(w_{+},u)_0\frac{dy}{y^n} = \frac{1}{2}
 \left[\frac{\|u\|_0^2}{y^{n-1}}\right]_a^b - \epsilon
 \int_a^b\|u\|_0^2\frac{dy}{y^n}.
$$
Therefore (\ref{eq:Imwyyn-1inte}) is computed as
\begin{equation}
  {\rm Im}\,\left[\frac{(w_{+},u)_0}{y^{n-1}}\right]_a^b =  k\,
 \left[\frac{\|u\|_0^2}{y^{n-1}}\right]_a^b - 2\epsilon k
 \int_a^b\|u\|_0^2\frac{dy}{y^n}
 - {\rm Im}\int_a^b(f,u)_0\frac{dy}{y^n},
 \nonumber
\end{equation}
which implies
 \begin{equation}
  {\rm Im}\,\left[\frac{(w_{+},u)_0}{y^{n-1}}\right]_a^b \leq k\,
 \left[\frac{\|u\|_0^2}{y^{n-1}}\right]_a^b 
 + C\|f\|_{\mathcal B}\|u\|_{{\mathcal B}^{\ast}}.
 \nonumber
\end{equation}
Note that for $z \not\in {\bf R}$, $w_{+}$ and $u$ are in $L^2((0,\infty);L^2(M_0);dy/y^n)$. Hence, there exists a sequence $b_1 < b_2 <\cdots \to \infty$ such that
$$
\frac{|(w_{+},u)(b_m)| + \|u(b_m)\|^2}{b_m^{n-1}} \to 0.
$$
For $w_+$, we take $a = y < b = b_m$  to have
$$
\frac{\|u(y)\|^2}{y^{n-1}} \leq C_k\left(\frac{\|w_+(y)\|^2}{y^{n-1}} 
+ \frac{|(w_{+},u)(b_m)| + \|u(b_m)\|^2}{b_m^{n-1}} 
+  \|f\|_{\mathcal B}\|u\|_{{\mathcal B}^{\ast}}\right).
$$
Letting $m \to \infty$, we see that
$$
 \frac{\|u(y)\|^2}{y^{n-1}} \leq C\left(\frac{\|w_+(y)\|^2}{y^{n-1}}
 + \|f\|_{\mathcal B}\|u\|_{{\mathcal B}^{\ast}}\right).$$
Dividing by $y$ and integrating from $1/R$ to $R$, we have
$$
 \frac{1}{\log R}\int_{1/R}^R\|u(y)\|^2 \frac{dy}{y^n}
 \leq \frac{C}{\log R}\int_{1/R}^R\|w_+(y)\|^2
 \frac{dy}{y^n} + 
  C\|f\|_{\mathcal B}\|u\|_{{\mathcal B}^{\ast}},
$$
which implies 
$$
\|u\|_{{\mathcal B}^{\ast}} \leq C\|w_+\|_{{\mathcal B}^{\ast}} + C\|f\|_{\mathcal B}\|u\|_{{\mathcal B}^{\ast}}.
$$
This, together with Lemma 4.7, yields 
$$
\|u\|_{{\mathcal B}^{\ast}} \leq C\|f\|_{\mathcal B}, \quad \forall z \in J_+.
$$
Similarly, we can prove the theorem for $z \in J_-$.
 \qed


\subsection{Radiation condition and uniqueness theorem}

The following theorem gives the fastest decay order of non-trivial solutions to the Helmholtz 
equation $(H_{free} - \lambda)u = 0$.


\begin{theorem} Let $\lambda > 0$.
If $u \in {\mathcal B}^{\ast}$ satisfies $(H_{free} - \lambda)u = 0$ for $0 < y < y_0$ with some 
 $y_0 > 0$, and
\begin{equation}
  \liminf_{R\to\infty}\frac{1}{\log R}\int_{1/R}^{y_0}\|u(y)\|^2_{L^2(M_0)}
  \frac{dy}{y^n} = 0,
 \nonumber
\end{equation}
then $u = 0$ for $0 < y < y_0$.
\end{theorem}

We stress that we have only to assume the equation $(H_{free} - \lambda)u = 0$ to be satisfied near $y = 0$. Theorem 4.9 and the unique continuation theorem yield the following corollary.

 
\begin{cor}
 $\ \ \ \sigma_p(H_{free})\cap\big(0,\infty\big) = \emptyset$.
\end{cor}

 The proof of Theorem 4.9 is reduced to the following result on the growth property of solutions to an abstract differential equation.

 Let $X$ be a Hilbert space and consider the following differential equation for an $X$-valued function $u(t)$: 
\begin{equation}
- u^{\prime\prime}(t) + B(t)u(t) + V(t)u(t) - Eu(t) = P(t)u(t), \quad t > 0,
\label{eq:Chap2Sect3Diffeq}
\end{equation}
$E > 0$ being a constant. We assume the following. 

\medskip
\noindent
{\it (B-1)   $B(t)$ is a non-negative self-adjoint operator valued function with domain $D(B(t)) = D \subset X$ independent of $t > 0$. For each $x \in D$, the map $(0,\infty) \ni t \to B(t)x \in X$ is $C^1$, and there exist constants $t_0 > 0$ and $\delta > 0$ such that
\begin{equation}
t\frac{dB(t)}{dt} + (1 + \delta)B(t) \leq 0, \quad \forall t  > t_0.
\label{eq:Chap2Sect3tdBtdt}
\end{equation}
(B-2) For any fixed $t$, $V(t)$ is bounded self-adjoint on $X$ and satisfies
\begin{equation}
V(t) \in C^1((0,\infty);{\bf B}(X)),
\label{eq:Chap2Sect3VtC1}
\end{equation}
\begin{equation}
\frac{1}{t}\|V(t)\| + 
\big\|\frac{dV(t)}{dt}\big\| \leq C(1 + t)^{-1-\epsilon},
\quad \forall t \geq 1,
\label{eq:Chap2Sect3normVt}
\end{equation}
for some constants $C, \epsilon > 0$.\\
\noindent
(B-3) For any fixed $t$, $P(t)$ is a closed (not necessarily self-adjoint) operator on $X$ with domain $D(P(t)) \supset D$ satisfying
\begin{equation}
P(t)^{\ast}P(t) \leq C(1 + t)^{-2-2\epsilon}\big(B(t) + 1\big).
\label{eq:Chap2Sect3B1tastB1t}
\end{equation}
Moreover,
$$
{\rm Re}\,P(t) := \frac{1}{2}\left(P(t) + P(t)^{\ast}\right)
$$
is a bounded operator on $X$ and satsifies}
\begin{equation}
\|{\rm Re}\,P(t)\| \leq C(1 + t)^{-1-\epsilon}, \quad \forall t > 0.
\label{eq:Chap2Sect3RealB1t}
\end{equation}


\begin{theorem}
Under the above assumptions (B-1), (B-2), (B-3), if
$$
\liminf_{t\to\infty}(\|u'(t)\|_X + \|u(t)\|_X) = 0
$$
holds, there exists $t_1 > 0$ such that $u(t) = 0$,  $\forall t > t_1$.
\end{theorem}

The proof of the above theorem is given in \cite{IsKu09}, Chap. 2, \S3.
Now we prove Theorem 4.9. If $u(y)$ satisfies
$$
- y^2(\partial_y^2 + \Delta_{free})u + (n-2)y\partial_yu - \frac{(n-1)^2}{4}u - \lambda u = 0, \quad 
0 < y < y_0<1,
$$
$v(t) = e^{(n-1)t/2}u(e^{-t})$ satisfies
$$
- \partial_t^2v(t) - e^{-2t}\Delta_{free}v(t) - \lambda v(t) = 0, \quad 
t > a = - \log y_0.
$$
Then the assumptions (B-1), (B-2), (B-3) are satisfied with $X = L^2(M_0)$,
$B(t) = - e^{-2t}\Delta_{free}$, $V(t) = P(t) = 0$, $\delta = 1$ and $t_0 > {\rm max}\,(1,a)$. 
By Lemma 4.2 (2),
$$
\lim_{T\to\infty}\frac{1}{T}\int_{t_0}^{T}\left(\|v(t\|^2 + \|v'(t)\|^2\right)dt = 0,
$$
which guarantees the assumption of Theorem 4.11. Hence $v(t) = 0$ for $t > t_1$ with some 
$t_1 > 0$ and, by the uniwue continuation, for $t >t_0$. This proves Theorem 4.9.

\medskip
We put
\begin{equation}
s(y) = \left\{
\begin{array}{cc}
1, & y < 1, \\
-1, & y > 1,
\end{array}
\right.
\label{S5s(y)}
\end{equation}
\begin{equation}
\widetilde\sigma_{\pm}(y,z) = \frac{n-1}{2} \mp i\sqrt{z}s(y).
\label{S5widetildesyz}
\end{equation}

We say that a solution $u \in {\mathcal B}^{\ast}$ of the equation $(H_{free} - \lambda)u = f$ satisfies the {\it outgoing radiation condition} (for $\widetilde\sigma_+$), or {\it incoming radiation condition} (for $\widetilde\sigma_-$), if
 the following (\ref{S5RadCond1/R1}) is fulfilled:
\begin{equation}
\lim_{R\to\infty}\frac{1}{\log R}\int_{1/R}^R
\|(D_y - \widetilde\sigma_{\pm}(y,\lambda))u(y)\|^2_{L^2(M_0)}\frac{dy}{y^n} = 0,
\label{S5RadCond1/R1}
\end{equation}

 
\begin{lemma}\label{S4RadCondUniqueness}
 Assume that $\lambda > 0$ and $u$ satisfies the equation $(H_{free} - \lambda)u = 0$, and the outgoing or incoming radiation condition. Then
 $u = 0$.
\end{lemma}

Proof. We assume that $u$ satisfies the outgoing radiation condition. We take 
\begin{equation}
\chi(t) = \left\{
\begin{array}{cc}
1 + t, &\quad - 1 \leq t < 0, \\
1 - t, &\quad 0 < t \leq  1, \\
0, & \quad |t| > 1,
\end{array}
\right.
\nonumber
\end{equation}
and put 
\begin{equation}
\rho(t) = \chi'(t) = \left\{
\begin{array}{cc}
1 , &\quad - 1 < t < 0, \\
-1, &\quad 0 < t < 1, \\
0, & \quad |t| > 1.
\end{array}
\right.
\label{S5rho=chi'}
\end{equation}
Taking $\chi_R(y) = \chi(\log y/\log R)$, we multiply the equation $(H_{free} - \lambda)u = 0$ by $\chi_R(y)\overline{u}$ and integrate over $M_{free}\times
(0,\infty)$ to obtain
\begin{equation}
\begin{split}
0 &= {\rm Im}\int_{0}^{\infty}(- D_y^2u + (n - 1)D_yu,\chi_Ru)_0\frac{dy}{y^n}  \\&=  {\rm Im}\, \frac{1}{\log R}\int_{0}^{\infty}\rho\big(\frac{\log y}{\log R}\big)\left(D_yu,u\right)_0\frac{dy}{y^{n}} \\
&= {\rm Im}\, \frac{1}{\log R}\int_{1/R}^{R}s(y)\left(D_yu,u\right)_0\frac{dy}{y^{n}}.
\end{split}
\nonumber
\end{equation}
Here we have used (\ref{S5s(y)}) and (\ref{S5rho=chi'}).
Letting $g = 
D_yu - \widetilde\sigma_+u$, we then have
\begin{equation} \label{3.09.2}
0 = {\rm Im}\,\frac{1}{\log R}\int_{1/R}^Rs(y)(g,u)_0\frac{dy}{y^n} + 
\frac{\sqrt{\lambda}}{\log R}\int_{1/R}^R\|u\|^2_0\frac{dy}{y^n}.
\end{equation}
By the Schwarz inequality, we have
$$
\left|\frac{1}{\log R}\int_{1/R}^Rs(y)(g,u)_0\frac{dy}{y^n}\right|
\leq \left(\frac{1}{\log R}\int_{1/R}^R\|g\|^2_0\frac{dy}{y^n}\right)^{1/2}\|u\|_{\mathcal B^{\ast}}
$$
The assumption of the lemma then implies that the first term in the right-hand side of (\ref{3.09.2})
tends to $0$ as $R \to \infty$. Thus,
$$
\lim_{R\to\infty}\frac{1}{\log R}\int_{1/R}^R\|u\|^2_0\frac{dy}{y^n} = 0.
$$
The lemma then follows from Theorem 4.9. \qed


\begin{lemma}\label{Radconds-1}
Let $u = R_{free}(z)f, \, z \in J_{\pm}$, and $s > 1/2$. Then there exists a constant $C=C_s(J_\pm) > 0$ 
such that the following inequality holds: \\
\begin{equation}
\|(D_y - \widetilde\sigma_{\pm}(y,z))u\|_{s-1} \leq C(\|u\|_{-s} + \|f\|_s), \quad \forall z \in J_{\pm}.
\nonumber
\end{equation}
\end{lemma}

Proof. We consider the outgoing case, $z \in J_+$. Take 
\begin{equation}
\varphi(y) = 
\left\{
\begin{array}{cc}
(- \log y)^{2s-1}, & \quad 0 < y < 1, \\
0,& \quad y > 1.
\end{array}
\right.
\nonumber
\end{equation} By Lemma 4.4, $- R_+ \leq - L_+$. Note that there exists a constant 
$0 < b=e^{(1-2s)/2} < 1$ such that $D_y\varphi + 2\varphi > 0$ if $y < b$. Therefore, 
letting $w_+ = (D_y - \widetilde\sigma_+(y,z))u$,
\begin{equation}
\begin{split}
& (2s-1)\int_a^b(-\log y)^{2s-2}\|w_+\|^2_0 \frac{dy}{y^n} \\
& \leq -\left[\frac{(-\log y)^{2s-1}(\|w_+\|^2_0 - \|D_xu\|^2_0)}{y^{n-1}}\right]_a^b 
+ 2\left|\int_a^b(-\log y)^{2s-1}(f,w_+)_0\frac{dy}{y^n}\right|.
\end{split}
\label{S5radcondcalculus1}
\end{equation}
By Lemma 4.2 (5), $w_+, D_x u \in L^{2,s}$. So, there exists a sequence 
$a_l \to 0,\, a_l \in \Bbb R_+,$ such that 
$$
\frac{(-\log y)^{2s-1}(\|w_+\|^2_0 - \|D_xu\|^2_0)}{y^{n-1}}\Big|_{y=a_l} \to 0. 
$$
Letting $a_l \to 0$ in (\ref{S5radcondcalculus1}), we then have
\begin{equation}
\begin{split}
&\|\theta(b-y) w_+\|^2_{s-1} \\
& \leq C_{s,\epsilon}(\|u(b)\|_0^2 + \|D_xu(b)\|^2_0 
 + \|D_yu(b)\|^2_0+ \|f\|^2_s) + \epsilon\|\theta(b-y)w_+\|^2_{s-1} \\
& \leq C_{s,\epsilon}(\|u\|_{-s}^2  + \|f\|^2_s) + \epsilon\|\theta(b-y)w_+\|_{s-1},
\end{split}
\nonumber
\end{equation}
where $\theta(\cdot)$ is the Heaviside function and
we have used that
$$
\|u(b)\|_0 + \|D_xu(b)\|_0 + \|D_yu(b)\|_0
\leq C_b(\|u\|_{\mathcal B^{\ast}} + \|f\|_{\mathcal B})
$$
by the standard a-pripri estimates for elliptic partial diffetential equations. This implies
\begin{equation}
\|\theta(b-y)w_+\|_{s-1} \leq C_s(\|u\|_{-s} + \|f\|_{s}).
\label{S5chiy<1w+estimate}
\end{equation}
Similarly, for $b'$ large enough, one can prove
\begin{equation}
\|\theta(y-b')w_+\|_{s-1} \leq C_s(\|u\|_{-s} + \|f\|_{s}).
\label{S5chiy>1w-estimate}
\end{equation}
At last, for $b<y <b'$, we use the standard elliptic estimates.
This completes the proof of the lemma. \qed


\begin{cor} \label{Cor:4.14}
Let $f \in L^{2,s}$ with $1/2 < s <1$. Then $R_{free}(z)f, \, z \in J_{+},$ satisfies the 
outgoing radiation condition, and $R_{free}(z)f, \,  z \in J_-,$ satisfies the 
incoming 
radiation condition.
\end{cor}

Proof.  Let $w_{\pm} = (D_y - \widetilde\sigma_{\pm}(y,z))u$ with $u = R_{free}(z)f$, $z \in J_{\pm}$. 
Then we have
\begin{equation}
\frac{1}{1+ \log R}\int_{1/R}^R\|w_+\|^2_0\frac{dy}{y^n} \leq 
(1 + \log R)^{1-2s}\|w_+\|^2_{s-1}.
\nonumber
\end{equation}
Letting $R \to \infty$ and taking into account of Lemma \ref{Radconds-1}, we obtain this corollary. \qed


\subsection{Limiting absorption principle}


\begin{theorem} \label{Th:4.15}
(1) For any $\lambda > 0$, $\lim_{\epsilon \to 0}
R_{free}(\lambda \pm i\epsilon) =: R_{free}(\lambda \pm i0)$ exists in the weak-$\ast$ sense, namely
\begin{equation}
\exists\lim_{\epsilon \to 0}(R_{free}(\lambda \pm i\epsilon)f,g) =: (R_{free}(\lambda \pm i0)f,g), \quad 
\forall f, g \in \mathcal B.
\nonumber
\end{equation} 
\noindent
(2) For any compact interval $I\subset (0,\infty)$ there exists a constant $C > 0$ such that
\begin{equation}
\|R_{free}(\lambda \pm i0)f\|_{{\mathcal B}^{\ast}} \leq C\|f\|_{\mathcal B}, \quad \forall \lambda \in I.
\nonumber
\end{equation}
(3) For any $\lambda > 0$ and $f \in \mathcal B$, $R_{free}(\lambda + i0)f$ satisfies the 
outgoing radiation condition, and $R_{free}(\lambda - i0)f$ satisfies the incoming radiation 
condition. \\
\noindent
(4) For any $f, g \in {\mathcal B}$, the map $(0,\infty) \ni \lambda \to (R_{free}(\lambda \pm i0)f,g)$ 
is continuous. 
\end{theorem}

The proof of the theorem is based upon the following lemmata.


\begin{lemma} \label{L:4.16}
Let $s > 1/2, \, f \in L^{2,s}$. \\
\noindent
(1) There exists a constant $C > 0$ such that
\begin{equation} \label{3.09.4}
\sup_{z \in J_{\pm}}\|R_{free}(z)f\|_{-s} \leq C\|f\|_{s}.
\end{equation}
(2) For any $\lambda > 0$ and $f \in L^{2,s}$, the strong limit $\lim_{\epsilon \to 0}R_{free}(\lambda \pm i\epsilon)f$ exists in $L^{2,-s}$. \\
\noindent
(3) $R_{free}(\lambda + i0)f$ satisfies the outgoing radiation condition, and 
$R_{free}(\lambda - i0)f$ satisfies the incoming radiation condition.\\
\noindent
(4) $\ R_{free}(\lambda \pm i0)f$ is an $L^{2,-s}$-valued continuous function of $\lambda > 0$. 
\end{lemma}

Proof. {\it Step 1}. Let $\lambda \in [a, b],\, f \in L^{2, s},\, s >1/2$. Assume $z_n \in J_+,\,
z_n \to \lambda$. Then there is a subsequence $\{z_{n'}\}$ of $\{z_n\}$ such that 
$R_{free}(z_{n'})f$ converges to some $u \in L^{2,-s}$. Indeed, take $1/2 < s' < s$ and let 
$\epsilon = (s - s')/2$. 
For the sequence $J_+ \ni z_n \to \lambda \in I$, we put $u_n = R_{free}(z_n)f$, and let 
$\chi_{R}(y)$ be such that $\chi_{R}(y) = 1$ for 
$y \not\in [1/R,R]$, $\chi_{R}(y) = 0$ for $y \in [1/R,R]$. By Theorem 4.8,
$$
\|\chi_{R}(y)u_n\|_{-s} \leq C(1 + \log R)^{-\epsilon}\|u_n\|_{-s'} \leq 
C(1 + \log R)^{-\epsilon}\|u_n\|_{\mathcal B^{\ast}}
\leq C(1 + \log R)^{-\epsilon},
$$
where $C$'s are constants independent of $n$ and $R$. This, together with 
Rellich's selection theorem in $H^2(M_0\times(1/R,R))$, implies that there exists a subsequence $\{u_{n'}\}$ of $\{u_n\}$ convergent to some $u \in L^{2,-s}$. 
Since $L^{2, s} \subset \mathcal B,\, \mathcal B^{\ast} \subset L^{2,-s}$, 
it follows from
Theorem 4.8,  that $u$ satisfies (\ref{3.09.4}).
Moreover, using again Theorem 4.8, we see that $u \in \mathcal B^{\ast}$ and
\begin{equation} \label{3.09.5}
\| u\|_{\mathcal B^{\ast}} \leq C \|f\|_{\mathcal B}.
\end{equation}
Indeed, from the definition of $\mathcal B^{\ast}$, for any $R>e$,
$$
\int_{1/R}^R \|u_{n'}\|_0^2 \frac{dy}{y^n} \leq C  \log R \|f\|_{\mathcal B}.
$$
Using the convergence of $u_{n'}$ to $u$ in $L^{2,-s}$, this inequality remains 
valid for $u$, implying
(\ref{3.09.5}). 

{\it Step 2}. By going to the limit, we see that
\begin{equation} \label{3.09.3}
(H_{free} - \lambda)u = f \quad {\rm on}\quad  \mathcal M_0.
\end{equation}
Since $z_{n'} \to \lambda $, this together with Lemma 4.2 (1), implies
\begin{equation} \label{3.09.6}
\|u_{n'}-u \|_{\mathcal B^{\ast}},\,\,\|D_y\left(u_{n'}-u \right)\|_{\mathcal B^{\ast}},
\,\,\|D_x\left(u_{n'}-u \right)\|_{\mathcal B^{\ast}} \to 0, \quad \hbox{as}\,\, n' \to \infty.
\end{equation}
Employing Corollary \ref{Cor:4.14} for $u_{n'}$, the equation (\ref{3.09.6}) implies that 
$u$ also satisfies
the outgoing radiation condition.

{\it Step 3}. The crucial point is that $u$ is independent of the choice of $z_{n'}$. 
Indeed, if there 
were another
limit, $u'$, then
$(H_{free}-\lambda) (u-u')=0$,
which, together with the outgoing radiation condition and the fact that 
$u-u' \in \mathcal B^{\ast} $, implies that $u=u'$, see Lemma 4.12. By this argument, every subsequence of $\{R_{free}(z_n)f\}$ contains a sub-subsequence which converges to one and the same limit. This shows that $\{R_{free}(z_n)\}$ itself converges without choosing subsequences.

Taking
$R_{free}(\lambda+i0) f= u$,
we see that
steps 2 and 3 yield statements (2) and (3) of the Lemma (the case of $R_{free}(\lambda -i0)$ can be
treated analogously).

{\it Step 4}. To prove (4),  assume the contrary and choose 
$f \in L^{2,s}$ and $\lambda_n \to \lambda$ such that, for some $\epsilon >0$,
$\|R_{free}(\lambda_n+i0) f- R_{free}(\lambda+i0)f \|_{-s} >\epsilon$.
However, there are $z_n \in J_+$ such that $\|R_{free}(\lambda_n+i0) f-R_{free}(z_n) f\|_{-s} < \epsilon/2$
and $|z_n-\lambda_n|< 1/n$. Since $R_{free}(z_n)f \to R_{free}(\lambda+i0)f$, we come to a contradiction.
 \qed

\bigskip
We extend Lemma 4.16  to $\mathcal B$ in the following way.


\begin{lemma} \label{L:4.18}
Let $f \in \mathcal B$. \\
\noindent
(1) For any $\lambda > 0$ and $f \in \mathcal B$, the weak$^*$ limit $\lim_{\epsilon \to 0}R_{free}(\lambda \pm i\epsilon)f$ exists in $\mathcal B^{\ast}$.
Moreover, there exists a constant $C > 0$ such that
\begin{equation} \label{4.09.1}
\sup_{z \in {\overline{J_{\pm}}}}\|R_{free}(z)f\|_{\mathcal B^{\ast}} \leq C\|f\|_{\mathcal B}.
\end{equation}

\noindent
For any 
 $s > 1/2$, \\
\noindent
(2) For any $\lambda > 0$ and $f \in \mathcal B$, the strong limit $\lim_{\epsilon \to 0}R_{free}(\lambda \pm i\epsilon)f$ exists in $L^{2,-s}$. \\
\noindent
(3) $R_{free}(\lambda + i0)f$ satisfies the outgoing radiation condition, and 
$R_{free}(\lambda - i0)f$ satisfies the incoming radiation condition.\\
\noindent
(4) $\ R_{free}(\lambda \pm i0)f$ is an $L^{2,-s}$-valued continuous function of $\lambda > 0$. 
\end{lemma}

Proof. Taking into account Theorem 4.8, it is sufficient to consider what happens when 
$z_n \to \lambda+i0$. First note that, for $f \in L^{2, s},\, s >1/2$, in the weak $\ast$-sense, 
\begin{equation} \label {4.09.2}
R_{free}(z_n) f \to R_{free}(\lambda+i0)f,\quad
\hbox{i.e.}\,\,
\left(R_{free}(z_n) f - R_{free}(\lambda+i0) f, g  \right) \to 0,
\end{equation}
for any $g \in \mathcal B$. For $g \in L^{2,s}$, (\ref{4.09.2}) follows from Lemma \ref{L:4.16} (2).
Taking into  account  of Theorem 4.8 and (\ref{3.09.5}), where, as we now know
$u= R_{free}(\lambda+i0) f$, and approximating   $g \in \mathcal B^{\ast}$ by $g_m \in L^{2,s}$,
we obtain (\ref{4.09.2}) for any $g \in \mathcal B$.

Let now $f \in \mathcal B$ and take $f_m \to f,\, f_m \in L^{2, s}.$ Then, for any $g \in \mathcal B$
and $z_n \to \lambda+i0$,
\begin{equation} \begin{split}
&\big|\left((R_{free}(z_n)- R_{free}(z_{n'}))f,\, g  \right)\big| \\
&\leq \big|\left((R_{free}(z_n)- R_{free}(z_{n'}))(f-f_m),\, g  \right)\big|+
\big|\left((R_{free}(z_n)- R_{free}(z_{n'})) f_m,\, g  \right)\big|.
\nonumber
\end{split}
\end{equation}
By Theorem 4.8, the first term in rhs is bounded by 
$C \|f-f_m  \|_{\mathcal B} \,\|g\|_{\mathcal B}$. As for the second term, for any $m$ it tends to $0$
by (\ref{4.09.2}). This implies that $\left(R_{free}(z_n)f,\, g  \right)$ is a Cauchy sequence proving
the existence of the weak $\ast$ limit $R_{free}(\lambda+i0) f \in \mathcal B^{\ast}$. Moreover, using
Theorem 4.8, we obtain estimate (\ref{4.09.1}).

To prove (2) and (4), we note that
$R_{free}(z) f= R_{free}(f_m) +R_{free}(z) (f-f_m)$,
where, as earlier, $f_m \in L^{2,s}$ approximate $f\in \mathcal B$. Since 
$\|u\|_{L^{2,-s}} \leq C \|u\|_{ \mathcal B^{\ast}}$, the claims follow from 
Lemma \ref{L:4.16} (2) and (4) and estimate (\ref{4.09.1}).

At last, using Lemma 4.2 (1) and Lemma \ref{L:4.16} (3), we see that $R_{free}(\lambda+i0) f$
satisfies the outgoing radiation condition. \qed

\bigskip
To complete the proof of Theorem 4.15, it remains to obtain statement (4). 
However, due to Lemma \ref{L:4.18} (4), this is true for $g \in L^{2,s},\, s>1/2$. 
Since $L^{2,s}$ is dense in $\mathcal B$, the case $g \in \mathcal B$ follows 
 from (\ref{4.09.1}). \qed

\bigskip
The following lemma is a consequence of the above proof and Lemma 4.12.


\begin{lemma}
 For any $f \in {\mathcal B}$ and $\lambda > 0$, $u = R_{free}(\lambda \pm i0)f$ satisfies the equation $(H_0 - \lambda)u = f$, 
and the radiation condition. 
Conversely, any solution $u \in {\mathcal B}^{\ast}$ of the above equation 
satisfying the radiation condition is unique and is given by $u = R_{free}(\lambda \pm i0)f$. 
 \end{lemma}


\subsection{Fourier transforms on the model space}
We construct the Fourier transform on $\mathcal M_{free}$. By projecting onto the eigenspace of $- \Delta_{free}$, $H_{free}$ becomes 
$$
H_{free}(\lambda_m) := y^2(- \partial_y^2 + \lambda_m) + (n-2) y\partial_y - \frac{(n-1)^2}{4},
$$
and this reduces the problem to 1-dimension.
We put
\begin{equation}
\hat f_m(y) = \int_{M_{free}}f(x,y)\overline{\varphi_m(x)}\,dV_h(x),
\label{S5hatfmy}
\end{equation}
\begin{equation}
\omega_{\pm}(k) = \frac{\pi}{(2k\sinh(k\pi))^{1/2}\Gamma(1 \mp ik)},
\label{Stomegapmk}
\end{equation}
where $\Gamma(z)$ is the gamma function, and
\begin{equation}
C_m^{(\pm)}(k) = 
\left\{
\begin{split}
\Big(\frac{\sqrt{\lambda_m}}{2}\Big)^{\mp ik} \quad (\lambda_m \neq 0), \\
\frac{\pm i}{k\omega_{\pm}(k)}\sqrt{\frac{\pi}{2}} \quad (\lambda_m = 0).
\end{split}
\right.
\label{C,pmk}
\end{equation}

Two ends of $\mathcal M_{free}$, the cusp and the regular infinity, give different contributions to the Fourier transforms. The part due to cusp is as follows:
\begin{equation}
\mathcal F_{c,free}^{(\pm)}(k)f = \frac{1}{\sqrt{|M_{free}|}}F_{free,0}^{(\mp)}(k)f,
\label{S5focfk}
\end{equation}
where $|M_{free}|$ is the volume of $M_{free}$, and
\begin{equation}
F^{(\pm)}_{free,0}(k)f = \frac{1}{\sqrt{2\pi}}\int_0^{\infty}y^{1\pm ik}\hat f_0(y)\frac{dy}{y^n}.
\label{S5Fpm0fk}
\end{equation}
The part due to regular end is as follows:
\begin{equation}
\mathcal F_{reg,free}^{(\pm)}(k)f = \sum_{m=0}^{\infty}
C_m^{(\pm)}(k)\varphi_m(x)F_{free,m}^{(\pm)}(k)f,
\label{S5foregpmfk}
\end{equation}
\begin{equation}
F_{free,m}^{(\pm)}(k) = \left\{
\begin{split}
& F_{free,m}(k)  \quad (\lambda_m \neq 0), \\
& F_{free}^{(\pm)}(k)  \quad (\lambda_m = 0),
\end{split}
\right.
\label{StFompm}
\end{equation}
where for $\lambda_m \neq 0$,
\begin{equation} \label{12.09.4}
F_{free,m}(k)f= \frac{(2k\sinh(k\pi))^{1/2}}{\pi}
\int_0^{\infty}y^{(n-1)/2} K_{ik}(\sqrt{\lambda_m}y)\hat f_m(y)\frac{dy}{y^n},
\end{equation}
$K_{\nu}(z)$ being the modified Bessel function.
We put
\begin{equation}
\big(\mathcal F_{c,free}^{(\pm)}f\big)(k) = \mathcal F_{c,free}^{(\pm)}(k)f, \quad 
\big(\mathcal F_{reg,free}^{(\pm)}f\big)(k) = \mathcal F_{reg,free}^{(\pm)}(k)f,
\end{equation}
\begin{equation}
\mathcal F_{free}^{(\pm)} = (\mathcal F_{c,free}^{(\pm)},\mathcal F_{reg,free}^{(\pm)}),
\label{S5mathcalF0pm}
\end{equation}
\begin{equation}
{\bf h}_{free} = {\bf C}\oplus L^2(M_{free}),
\label{S5bfh}
\end{equation}
\begin{equation}
\widehat{\mathcal H}_{free} = L^2((0,\infty) ; {\bf h}_{free} ; dk),
\label{S5widehatH}
\end{equation}
where (\ref{S5widehatH}) is the set of all ${\bf h}_{free}$-valued $L^2$-functions on $(0,\infty)$ with respect to the measure $dk$.
Then we can obtain the following theorem exactly in the same way as 
\cite{IsKu09}, Chap. 3, Theorem 2.5.


\begin{theorem}
$\mathcal F_{free}^{(\pm)}$, defined on $C_0^{\infty}(\mathcal M_{free})$, is uniquely extended to a unitary operator from $L^2(\mathcal M_{free})$ to $\widehat{\mathcal H}_{free}$. Moreover, for $f \in D(H_{free})$,
$$
(\mathcal F_{free}^{(\pm)}H_{free}f)(k) = k^2(\mathcal F_{free}^{(\pm)}f)(k).
$$
\end{theorem}


\section{Laplace-Beltrami operators on orbifolds}

In the case of the orbifold with asymptotically hyperbolic ends, $\mathcal M$ of form
(\ref{S1MKM1MN}), with ${\mathcal M}_j$ satisfying {\bf (A-1)}--{\bf (A-3)}, we use a finite uniformising cover which, at every end ${\mathcal M}_j$ is the product form described in the beginning of \S 4. We denote the corresponding partition of unity by $\chi_j(X),
\, j=1, \dots, m$.


\subsection{Spectral properties of the Laplace-Beltrami operator}
Let $H$ be the shifted Laplacian on $\mathcal M$, 
\begin{equation}
H = - \Delta_g - \frac{(n-1)^2}{4}.
\label{S4H}
\end{equation}
Under the assumption {\bf(A-3)}, it may be written, in any ${\mathcal M}_j$, as a perturbation of the "unperturbed" operator described in \S 4
\begin{equation}
H_{free(j)} = - y^2(\partial_y^2 + \Delta_{j}) + (n-2) y\partial y - \frac{(n-1)^2}{4},
\label{S4HinMi}
\end{equation}
where $\Delta_{j}= \Delta_{h_j}$.


\begin{theorem}
Under the assumptions {\bf (A-1)}--{\bf (A-3)}, $-\Delta_g\big|_{C_0^{\infty}(\mathcal M)}$ is essentially self-adjoint.
\end{theorem}

Sketch of proof. Using the assumptions {\bf (A-1)}--{\bf (A-3)} and following the proof of	Theorem 3.2 in \cite{IsKu09}, we show that if
$$
(-\Delta_g-z) u=f, \quad u, f \in L^2(\mathcal M),\quad \hbox{Im}\, z \neq 0,
$$
then
\begin{equation} 
\|u\|_{H^2(\mathcal M)} \leq C (1+|z|) (\|u\|+\|f\|).
\end{equation}
The only difference from the proof in \cite{IsKu09} is that, when dealing with the orbifold ends
we decompose $u$ as
$$
u= \sum_{j=1}^m \chi_j u,
$$
where $\chi_j$ is the special partition of unity subordinate to the orbifold structure
which is introduced in the beginning of this section.

Let $A = - \Delta_g\big|_{C_0^{\infty}(\mathcal M)}$, and  show that 
$N(A^{\ast} \pm i)) = \{0\}$. Suppose that, for some $u \in L^2(\mathcal M)$,
$(A^{\ast} + i) u = 0$. Taking $u_R=\chi_R u$, where $\chi_R$ equals to $1$ in $\mathcal K$ and
in $M_j \times (1, R),\, j=1, \dots N+N',$ and $0$ in $M_j \times (R+1, \infty)$ and looking at
$\hbox{Im}\left( (A^{\ast} + i) u_R, u_R \right)$, we get the result.
\qed


\begin{theorem} \label{spectrum}
(1) $\sigma_e(H) = [0,\infty)$.\\
\noindent
(2) $\sigma_{d}(H) \subset (-\infty,0)$.
\end{theorem}

Proof.  Similarly to the proof of Theorem 3.2 (2) in  \cite{IsKu09}, we
can represent the resolvent $R(z)$ of $H$ in the form
\begin{equation} \label{12.09.2}
R(z)= \sum_{j=1}^{N+N'} { \chi}_j R_{free(j)}(z) {\widetilde \chi}_j
+R(z) (\chi_0-A(z)).
\end{equation}
	Here the partition of unity $\chi_j$ is different from the one described in the beginning
of the section. Namely, we take 
\begin{equation} \label{12.09.3}
\chi_j=1\,\, \hbox{in}\,\,M_j \times (2, \infty),\quad \chi_j=0\,\, \hbox{outside}\,\,
M_j \times (3/2, \infty),  \quad \chi_0=1 -\sum_{j=1}^m \chi_j.
\end{equation}
 We also take
${\widetilde \chi}_j=1$ on $\hbox{supp}(\chi_j)$ and $0$ outside $M_j \times (1, \infty)$.
As for $R_{free(j)}$, this is the resolvent of $H_{free(j)}$, and
\begin{equation} \label{12.09.6}
\begin{split}
A(z) & =\sum_{j=1}^{N+N'} A_j(z) {\widetilde \chi}_j, \\
A_j(z) & =[H, \chi_j] R_{free(j)}(z) +\chi_j (H-H_{free(j)}) {\widetilde \chi}_j  R_{free(j)}(z).
\end{split} 
\end{equation}
Then the proof of the theorem follows the same arguments as that 
of Theorem 3.2 (2) in \cite{IsKu09}\qed

\begin{remark}
When $N'>0$, $\sigma_p(H) \cap (0, \infty)= \emptyset.$ This can be proven
as in Theorem 3.5 (1) in  \cite{IsKu09}.
\end{remark}

Using partition  $\{\chi_j\}$ and  defining, in a natural way, the norms $\|\cdot\|_{\mathcal B}$, $\|\cdot\|_{\mathcal B^{\ast}}$ on $\mathcal M_i$, we  put
$$
\|f\|_{\mathcal B} = \|\chi_0f\|_{L^2(\mathcal M)} + \sum_{j=1}^N\|\chi_j f\|_{\mathcal B},
$$
$$
\|u\|_{\mathcal B^{\ast}} = \|\chi_0u\|_{L^2(\mathcal M)} + \sum_{j=1}^N\|\chi_j u\|_{\mathcal B^{\ast}},
$$
which define the Besov-type spaces $\mathcal B$ and $\mathcal B^{\ast}$ on $\mathcal M$. 

Once we have established the resolvent estimates for the model space, we can follow the arguments for the spectral properties of the perturbed operator $H$ in the same way as in \cite{IsKu09} by using the partition of unity as above. Henceforth, we state only the results quoting the corresponding theorems in \cite{IsKu09}.

 
\begin{theorem} \label{Th:6.2}
For $\lambda \in \sigma_e(H)\setminus\sigma_p(H)$, there exists a limit
$$
\lim_{\epsilon \to 0}R(\lambda \pm i\epsilon) \equiv 
R(\lambda \pm i0) \in {\bf B}({\mathcal B};{\mathcal B}^{\ast})
$$
in the weak $\ast$-sense. Moreover for any compact interval $I \subset \sigma_e(H)\setminus\sigma_p(H)$ there exists a constant $C > 0$ such that
\begin{equation}
\|R(\lambda \pm i0)f\|_{{\mathcal B}^{\ast}} \leq C\|f\|_{\mathcal B},
\quad \lambda \in I.
\nonumber
\end{equation} 
 For $f, g \in {\mathcal B}$,  $(R(\lambda \pm i0)f,g)$ is continuous with respect to $\lambda > 0$.
\end{theorem}

The proof is similar to that of Theorem 3.8 in \cite{IsKu09}.

\medskip
Next we introduce the {\it radiation conditions} on $\mathcal M$. Let 
$$
\sigma_{\pm}(\lambda) = \frac{n-1}{2} \mp i\sqrt{\lambda}, \quad \lambda > 0.
$$
We say that a solution $u \in \mathcal B^{\ast}$ of the equation $(H - \lambda)u = f \in \mathcal B$ satisfies the {\it outgoing radiation condition}, or $u$ is {\it outgoing}, if
\begin{eqnarray*}
& \displaystyle \lim_{R\to\infty}\frac{1}{\log R}\int_{2}^{R}\|\big(y\partial_y - \sigma_-(\lambda)\big)u(\cdot,y)\|^2_{L^2(M_i)}\frac{dy}{y^n} = 0, \quad (i = 1, \cdots, N), \\
&\displaystyle \lim_{R\to\infty}\frac{1}{\log R}\int_{1/R}^{1/2}\|\big(y\partial_y - \sigma_+(\lambda)\big)u(\cdot,y)\|^2_{L^2(M_i)}\frac{dy}{y^n} = 0, \quad (i = N+1, \cdots, N+N'),
\end{eqnarray*}
and similarly for the {\it incoming radiation condition}.

The following theorem is analogous to Theorem 3.7 in \cite{IsKu09}.


\begin{theorem}\label{S5uniqueness}
Let $\lambda \in (0,\infty)\setminus\sigma_p(H)$, and suppose $u \in \mathcal B^{\ast}$ satisfies $(H - \lambda) u = 0$ and the radiation condition. Then : \\
\noindent
(1) If one of $\mathcal M_j$ has a regular infinity, then $u=0$. \\
\noindent
(2) If all $\mathcal M_j$ has a cusp, then $u \in L^{2,s}, \ \forall  s > 0$.
.
\end{theorem}


\subsection{Fourier transforms associated with $H$} 

\subsubsection{Definition of ${\mathcal F}_{free(j)}^{(\pm)}(k)$}

$\ $ \\
\noindent
(i) For $1 \leq j \leq N$ (the case of cusp), let $\mathcal F_{c,free(j)}^{(\pm)}(k)$ be defined by  (\ref{S5focfk}) with ${ M}_j$ in place of ${ M}_{free}$.

\noindent 
(ii) For $N+1 \leq j \leq N+N'$ (the case of regular infinity), let $\mathcal F^{(\pm)}_{reg,free(j)}(k)$ be defined by (\ref{S5foregpmfk}) with $M_{free}$ replaced by $M_j$, and $\lambda_m, \varphi_m(x)$ by the eigenvalues and complete orthonormal system of eigenvectors of $- \Delta_j$ : 
$$
\lambda_{j,0} < \lambda_{j,1} \leq \cdots ; \quad
\varphi_{j,0}(x), \ \varphi_{j,1}(x), \ \cdots.
$$

\subsubsection{ Definition of ${\mathcal F}^{(\pm)}(k)$.} 
 The Fourier transform associated with $H$ is now defined by
\begin{equation}
{\mathcal F}^{(\pm)}(k) = 
\big({\mathcal F}_{1}^{(\pm)}(k), \cdots, {\mathcal F}_{N+N'}^{(\pm)}(k)\big),
\label{S3Fpmk}
\end{equation}
where,
for $1 \leq j\leq N+N'$,
\begin{equation}
\begin{split}
{\mathcal F}_j^{(\pm)}(k) =  
{\mathcal F}_{free(j)}^{(\pm)}(k)Q_j(k^2 \pm i0),
\end{split}
\label{eq:Fjplusminuskcusp}
\end{equation}
where  denoting
\begin{equation}
\widetilde V_j = H - H_{free(j)} \quad {\rm on} \quad \mathcal M_j,
\label{S3tildeVj}
\end{equation}
we put
\begin{equation}
Q_j(z) = \chi_j + \Big([H_{free(j)},\chi_j] - \chi_j\widetilde V_j\Big)R(z).
\label{S3DefineQ(z)}
\end{equation}

\medskip
For functions $f, g \in {\mathcal B}^{\ast}$ on ${\mathcal M}$, by
$f \simeq g$
we mean that on each end $\mathcal M_j$
\begin{equation}
\begin{split}
\lim_{R\to\infty}\frac{1}{\log R}\int_{1 < y < R}\|f(y) - g(y)\|^2_{L^2(M_j)}
\frac{dy}{y^n} = 0, \quad 1 \leq j\leq N, \\
\lim_{R\to\infty}\frac{1}{\log R}\int_{1/R < y < 1}\|f(y) - g(y)\|^2_{L^2(M_j)}
\frac{dy}{y^n} = 0, \quad N+1 \leq j \leq N+N'.
\end{split}
\nonumber
\end{equation}


\begin{theorem} Let $f \in {\mathcal B}$, 
$k^2 \in \sigma_e(H)\setminus\sigma_p(H)$. Then,
\begin{equation}
\begin{split}
R(k^2 \pm i0) f \simeq&\ \omega_{\pm}^{(c)}(k)\sum_{j=1}^N
\chi_j y^{(n-1)/2\pm ik} {\mathcal F}^{(\pm)}_j(k)f\\
&+ \omega_{\pm}(k)\sum_{j=N+1}^{N+N'}\chi_jy^{(n-1)/2 \mp ik}\mathcal F_j^{(\pm)}(k)f.
\end{split}
\nonumber
\end{equation}
\end{theorem}
The proof is similar to Theorem 3.10 in \cite{IsKu09}.

We put
\begin{equation}
{\bf h}_{\infty} = \left(\oplus_{j=1}^N{\bf C}\right) \oplus 
\left(\oplus_{j=N+ 1}^{N+N'}L^2(M_j)\right),
\label{eq:Chap3Sect2hinfty}
\end{equation}
and, for $\varphi, \psi \in {\bf h}_{\infty}$, we define the inner product by
\begin{equation}
(\varphi,\psi)_{{\bf h}_{\infty}} = \sum_{j=1}^N\varphi_j\overline{\psi_j}+
\sum_{j=N+1}^{N+N'}(\varphi_j,\psi_j)_{L^2(M_j)}.
\nonumber
\end{equation}
We put
\begin{equation}
\widehat{\mathcal H} = L^2((0,\infty);{\bf h}_{\infty};dk).
\nonumber
\end{equation}


\begin{theorem}
We define $\big({\mathcal F}^{(\pm)}f\big)(k) = {\mathcal F}^{(\pm)}(k)f$ for $f \in {\mathcal B}$. Then  ${\mathcal F}^{(\pm)}$ is uniquely extended to a bounded operator from $L^2({\mathcal M})$ to $\widehat{\mathcal H}$ with the following properties. \\
\noindent
(1) $\ {\rm Ran}\,\,{\mathcal F}^{(\pm)} = \widehat{\mathcal H}$. \\
\noindent
(2) $\ \|f\| = \|{\mathcal F}^{(\pm)}f\|$ for $f \in {\mathcal H}_{ac}(H)$. \\
\noindent
(3) $\ {\mathcal F}^{(\pm)}f = 0$ for $f \in {\mathcal H}_p(H)$.  \\
\noindent
(4) $\ 
\left({\mathcal F}^{(\pm)}Hf\right)(k) = 
k^2\left({\mathcal F}^{(\pm)}f\right)(k)$ for 
$f \in {\it D}(H)$. \\
\noindent
(5)$\ {\mathcal F}^{(\pm)}(k)^{\ast} \in {\bf B}({\bf h}_{\infty};{\mathcal B}^{\ast})$ and
$(H - k^2){\mathcal F}^{(\pm)}(k)^{\ast} = 0$  for $k^2 \in (0,\infty)\setminus\sigma_p(H)$. \\
\noindent
(6) For $f \in {\mathcal H}_{ac}(H)$, the inversion formula holds:
\begin{eqnarray*}
f = \left({\mathcal F}^{(\pm)}\right)^{\ast}{\mathcal F}^{(\pm)}f 
= \sum_{i=1}^{N+N'}\int_0^{\infty}{\mathcal F}^{(\pm)}_i(k)^{\ast}
\left({\mathcal F}^{(\pm)}_if\right)(k)dk.
\end{eqnarray*}
\end{theorem}
For the proof, see that of Theorem 3.12 in \cite{IsKu09}.


\subsection{$S$ matrix} 
We shall introduce the S-matrix in terms of the asymptotic expansion of solutions to the Helmholtz equation.


\begin{theorem}
If $k^2 \not\in (0,\infty)\sigma_p(H)$, we have 
\begin{equation}
{\mathcal F}^{(\pm)}(k){\mathcal B} = {\bf h}_{\infty},
\nonumber
\end{equation}
\begin{equation}
\{u \in {\mathcal B}^{\ast}\,;\,(H - k^2)u = 0\} = 
{\mathcal F}^{(\pm)}(k)^{\ast}{\bf h}_{\infty},
\nonumber
\end{equation}
cf. Theorem 3.13 in \cite{IsKu09}.
\end{theorem}

We derive an asymptotic expansion of solutions to the Helmholtz equation. Let $V_j$ be the differential operator defined by
\begin{equation}
 V_j = [H_{free(j)},\chi_j] + \widetilde V^{\ast}\chi_j \quad 
 (1 \leq j \leq N+N').
 \nonumber
\end{equation}
For $N+1 \leq j \leq N+N', \ 1 \leq \ell \leq N+N'$, we define
\begin{equation}
 \widehat S_{j\ell}(k) = \delta_{j\ell}J_j(k) - \frac{\pi i}{k}
 {\mathcal F}_j^{(+)}(k)V_{\ell}^{\ast}\left({\mathcal F}_{free(\ell)}^{(-)}(k)\right)^{\ast},
\nonumber
\end{equation}
\begin{equation}
J_j(k)\psi = \sum_{m\geq0}
\left(\frac{\sqrt{\lambda_{j,m}}}{2}\right)^{-2ik}
\varphi_{j,m}(x)
\widehat\psi_m \quad
(N+1 \leq j \leq N+N').
\nonumber
\end{equation}
For $1 \leq j \leq N, 1 \leq \ell \leq N+N'$, we define
\begin{equation}
\widehat S_{j\ell}(k) = - \frac{\pi i}{k}{\mathcal F}_j^{(+)}(k)V_{\ell}^{\ast}\left({\mathcal F}_{free(\ell)}^{(-)}(k)\right)^{\ast}.
\nonumber
\end{equation}


\begin{theorem}
For $\psi = (\psi_1,\cdots,\psi_N) \in {\bf h}_{\infty}$
\begin{equation}
 \begin{split}
   \left({\mathcal F}^{(-)}(k)\right)^{\ast}\psi
  \ \simeq &   \frac{ik}{\pi} \omega_-^{(c)}(k)\sum_{j=1}^N
   \chi_j y^{1-ik}\widehat \psi_{j0} \\
+&
  \frac{ik}{\pi} \omega_-(k)\sum_{j=N+1}^{N+N'}\chi_j 
  y^{(n-1)/2+ik}\psi_j    \\
    - &\frac{ik}{\pi} \omega_+^{(c)}(k)\sum_{j=1}^N\sum_{\ell=1}^{N+N'}
  \chi_j y^{(n-1)/2+ik}
  \widehat S_{j\ell}(k)\psi_{\ell} \\
   - &\frac{ik}{\pi} \omega_+(k)\sum_{j=N+1}^{N+N'}\sum_{\ell=1}^{N+N'}
  \chi_j y^{(n-1)/2-ik}
  \widehat S_{j\ell}(k)\psi_{\ell}.
   \end{split}
   \nonumber
\end{equation}
\end{theorem}
Cf. Theorem 3.14 in \cite{IsKu09}.

\medskip
We define an operator-valued $(N+N') \times (N+N')$ matrix $\widehat S(k)$ by
\begin{equation}
\widehat S(k) = \Big(\widehat S_{j\ell}(k)\Big)_{j, \ell=1}^{N+N'},
\nonumber
\end{equation}
and call it  $S$-matrix.  


\begin{theorem}
(1)  For any $u \in {\mathcal B}^{\ast}$ satisfying $(H - k^2)u = 0$, there exists a unique $\psi^{(\pm)} \in {\bf h}_{\infty}$ such that
\begin{equation}
 \begin{split}
    u &\simeq \omega_-(k)\sum_{j=N+1}^{N+N'}\chi_j 
  y^{(n-1)/2+ik}
  \psi_j^{(-)}  + \omega_-^{(c)}(k)\sum_{j=1}^N
   \chi_j y^{(n-1)/2-ik}\widehat\psi_{j0}^{(-)} \\
    &\  \ -  \omega_+(k)\sum_{j=N+1}^{N+N'}
  \chi_j y^{(n-1)/2-ik}
  \psi_j^{(+)} -  \omega_+^{(c)}(k)\sum_{j=N+1}^{N+N'}
  \chi_j y^{(n-1)/2+ik}
  \widehat\psi_{j0}^{(+)}.
   \end{split}
   \nonumber
\end{equation}
(2) For any $\psi^{(-)} \in {\bf h}_{\infty}$, there exists a unique $\psi^{(+)} \in {\bf h}_{\infty}$ and $u \in \mathcal B^{\ast}$ satisfying $(H - k^2)u = 0$, for which the expansion (1) holds. Moreover
\begin{equation}
\psi^{(+)} = \widehat S(k)\psi^{(-)}.
\nonumber
\end{equation}
(3) $\widehat S(k)$ is unitary on ${\bf h}_{\infty}$.
\end{theorem}
Cf. Theorem 3.15 in  \cite{IsKu09}.


\section{Generalized S-matrix}


\subsection{Exponentially growing solutions}

\begin{definition}\label{DefSequentialSpaces}
We introduce the sequential spaces $\ell^{2,\pm \infty}$ by
$$
\ell^{2,\infty} \ni a = (a_m)_{m\in{\bf Z}_+} \Longleftrightarrow 
\sum_{m\in{\bf Z}_+}|a_m|^2\rho^{m} < \infty, \quad \forall \rho > 1,
$$
$$
\ell^{2,-\infty} \ni b = (b_m)_{m\in{\bf Z}_+} \Longleftrightarrow 
\sum_{m\in{\bf Z}_+}|b_m|^2\rho^{-m} < \infty, \quad \exists \rho > 1,
$$
\end{definition}
cf. Definition 4.1 in \S3, \cite{IsKu09}.

Let $0 \neq k \in {\bf R}$. Suppose $u(x,y) \in C^{\infty}(M_j \times (1, \infty))$ 
 and satisfies the 
the equation
\begin{equation}
- y^2\big(\partial_y^2 + \Delta_j\big)u - \frac{(n-1)^2}{4}u = k^2 y, \quad 
y > 1.
\label{S4freeeq}
\end{equation}
 Expanding $u$ into a Fourier series
$$
u(x,y) = \sum_{m\in{\bf Z}_+} u_m(y) \varphi_{j,m}(x),
$$
we have
$$
y^2\big(- \partial_y^2 + \lambda_{j,m}\big)u_m - \frac{(n-1)^2}{4}u_m = k^2 u_m, \quad y > 1.
$$
Then $u_m$ can be written as
\begin{equation}
u_m(y) = \left\{
\begin{split}
& a_m\, y^{(n-1)/2} I_{-ik}({\sqrt{\lambda_{j,m}}}y) + b_{m}\, y^{(n-1)/2} K_{ik}({\sqrt{\lambda_{j,m}}} y), 
\quad (m\neq 0), \\
& a_0\, y^{(n-1)/2 - ik} + b_0\,y^{(n-1)/2 + ik}, \quad (m = 0).
\end{split}
\right.
\label{S4anbn}
\end{equation}
Here let us note that $K_{-\nu}(z) = K_{\nu}(z)$, and $K_{\nu}(z)$, $I_{\nu}(z)$ are linearly independent solutions to the equation
$$
\frac{d^2w}{dz^2} + \frac{1}{z}\frac{dw}{dz}-\left(1 + \frac{\nu^2}{z^2}\right)w = 0.
$$


\begin{lemma}
Let $a = (a_m)_{m\in{\bf Z}_+}$, $b = (b_m)_{m\in{\bf Z}_+}$ be defined by (\ref{S4anbn}).
If $a  \in l^{2,\infty}$, then $b \in l^{2,-\infty}$.
\end{lemma}
Proof. 
Recall the asymptotic expansion of modified Bessel functions
\begin{equation}
\left\{
\begin{split}
& I_{\nu}(z) \sim \frac{1}{\sqrt{2\pi z}}e^{z}, \quad z \to \infty, \\
& K_{\nu}(z) \sim \sqrt{\frac{\pi}{2z}}e^{-z}, \quad z \to \infty.
\end{split}
\right.
\label{S4BesselAsymp}
\end{equation}
Since $a \in \ell^{2,\infty}$, we have 
$\sum_{m\geq0}|a_m|^2\big|I_{-ik}({\sqrt{\lambda_{j,m}}} y)\big|^2 < \infty$ for any $y > 0$. 
By Parseval's formula,
$$
y^{1-n}\|u(\cdot,y)\|^2_{L^2(M_j)} = \sum_{m\geq 0}\big|a_mI_{-ik}
({\sqrt{\lambda_{j,m}}} y) + 
b_mK_{ik}({\sqrt{\lambda_{j,m}}} y)\big|^2 + |a_0y^{-ik} + b_0y^{ik}|^2.
$$
We then have $\sum_{m\geq0}|b_m|^2\big|K_{ik}({\sqrt{\lambda_{j,m}}} y)\big|^2 < \infty$, hence 
$b \in \ell^{2,-\infty}$. \qed

\medskip
We introduce the {\it spaces of generalized scattering data at infinity} :
\begin{equation}
{\bf A}_{\pm\infty} = \left({\mathop\oplus_{j=1}^N}\ell^{2,\pm\infty}\right)
\oplus \left({\mathop\oplus_{j=N+1}^{N+N'}}L^2(M_j)\right).
\label{S4bfApm}
\end{equation}

\bigskip
We use the following notation. For
\begin{equation}
\psi^{(-)} = (a_{1},\cdots,a_N, \psi_{N+1}^{(-)},\,\cdots,\psi^{(-)}_{N+N'}) \in {\bf A}_{\infty},
\label{S4psi-}
\end{equation}
\begin{equation}
\psi^{(+)} = (b_{1},\cdots,b_N,\,\psi_{N+1}^{(+)},\cdots,\psi^{(+)}_{N+N'},) \in {\bf A}_{-\infty},\quad 1 \leq j \leq N,
\label{S4psi+}
\end{equation}
let 
\begin{equation}
u_j^{(-)} = 
\left\{
\begin{split}
&\omega_-^{(c)}(k)\left(a_{j,0}\,y^{(n-1)/2-ik} + \sum_{m=1}^\infty
a_{j,m}\, \varphi_{j,m}(x) y^{(n-1)/2}I_{-ik}({\sqrt{\lambda_{j,m}}} y)\right), \quad 1 \leq j \leq N, \\
&\omega_-(k)\,y^{(n-1)/2+ik}\psi_j^{(-)}(x),
\end{split}
\right.
\label{S4uj-}
\end{equation}
\begin{equation}
u_j^{(+)} = 
\left\{
\begin{split}
&\omega_+^{(c)}(k)\left(b_{j,0}\,y^{(n-1)/2+ik} + \sum_{m=1}^\infty
b_{j,m}\,\varphi_{j,m}(x) y^{(n-1)/2}K_{ik}({\sqrt{\lambda_{j,m}}} y)\right), \quad 1 \leq j \leq N,\\
&\omega_+(k)\, y^{(n-1)/2-ik}\psi_j^{(+)}(x), \quad N+1 \leq j \leq N+N',
\end{split}
\right.
\label{S4uj+}
\end{equation}
where $a_{j,m}, b_{j,m}$ are the $m$-th components of $a_j \in \ell^{2,\infty}, b_j \in \ell^{2,-\infty}$.


\begin{lemma}\label{S6LongLemma}
Take $k > 0$ such that $k^2 \not \in \sigma_p(H)$, and let $\psi^{(-)}$, $u_j^{(-)}$ be as in (\ref{S4psi-}), (\ref{S4uj-}).
Then there exists a unique solution $u$ such that 
$$
(H - k^2) u = 0, \quad u - \sum_{j=1}^N {\widehat \chi}_ju_j^{(-)} 
\ {\rm is} \ {\rm in} \   \mathcal B^{\ast},\ {\rm and} \ {\rm outging},
$$ 
where ${\widehat \chi}_j(y)=1$ for $y >3$ and $0$ for $y <2$.
  For this $u$, there exists $\psi^{(+)} \in {\bf A}_{\infty}$ and $u_j^{(+)}$ of the form (\ref{S4uj+}) such that \\
\noindent 
(1) 
 For $j = 1,\cdots,N$, there exists $y_0 > 0$ such that in $\mathcal M_j$,
\begin{equation}
u = u_j^{(-)} - u_j^{(+)}, \quad {\rm if} \quad y > y_0.
\label{S4uu-h+}
\end{equation}
(2)
For $j = N+1,\cdots, N+N',$
\begin{equation}
u - u_j^{(-)} \simeq - u_j^{(+)}, \quad {\rm in} \quad \mathcal M_j.
\label{S4uminusu(-)}
\end{equation}
Explicitly, $b_j$ is given by
\begin{equation}
b_{j,0} = \mathcal F_j^{(+)}(k)f,
\label{S4bj0}
\end{equation}
\begin{equation}
f = (H-k^2)u^{(-)}, \quad u^{(-)} = \sum_{j=1}^{N+N'}{\widehat \chi}_ju_j^{(-)}.
\label{S4f}
\end{equation}
\begin{equation}
b_{j,m} = \int_0^{\infty}y^{(n-1)/2}I_{-ik}({\sqrt{\lambda_{j,m}}} y)f_{j,m}(y)\frac{dy}{y^n}, 
\quad m \geq 1.
\label{S4bjn}
\end{equation}
\begin{equation}
f_{j,m} = \langle f_j, \varphi_{j,m}\rangle_{L^2(M_j)}, \quad
f_j = \chi_jf + [H_{free(j)},\chi_j]R(k^2 + i0)f.
\end{equation}
\end{lemma}

Proof. The uniqueness follows from Theorem \ref{S5uniqueness}. To prove the existence,
we put
\begin{equation}
u = u^{(-)} - R(k^2 + i0)f.
\label{S4u}
\end{equation}
 By Theorem 5.6, we then have
\begin{equation}
\begin{split}
R(k^2 + i0)f \simeq & \ \omega_+^{(c)}(k)\sum_{j=1}^N\chi_jy^{(n-1)/2+ik}\mathcal F_j^{(+)}(k)f \\
& + \omega_+(k)\sum_{j=N+1}^{N+N'}\chi_jy^{(n-1)/2-ik}\mathcal F_j^{(+)}(k)f .
\end{split}
\nonumber
\end{equation}
Letting $b_{j,0} = \mathcal F_j^{(+)}(k)f$, $(j = 1,\cdots,N)$ and
$\psi_j^{(+)} = \mathcal F_j^{(+)}(k)f$ $(j=N+1,\cdots, N+N')$,   
we  prove the existence of $u$ and 
$\psi_j^{(+)}$, $b_{j,0}$.

Direct computation shows that
$$
(H_{free(j)} - \lambda)\chi_jR(\lambda \pm i0) = \chi_j + [H_{free(j)},\chi_j]R(\lambda \pm i0),
$$
In fact, take $\widetilde\chi_j \in C^{\infty}(\mathcal M)$ such that ${\rm supp}\, \widetilde\chi_j \subset \mathcal M_j$ and $\widetilde \chi_j=1$ on ${\rm supp}\,\chi_j$. Then, since $ds^2 = (dy)^2 + h_j(\widetilde x,d\widetilde x)$ on $\mathcal M_j$, the left-hand side is equal to
\begin{equation}
\begin{split}
&\chi_j(H_{free(j)}-\lambda)\widetilde\chi_jR(\lambda \pm i0)  + 
[H_{free(j)},\chi_j]\widetilde\chi_jR(\lambda \pm i0) \\
=& \chi_j(H-\lambda)\widetilde\chi_jR(\lambda \pm i0) +  [H_{free(j)},\chi_j]\widetilde\chi_jR(\lambda \pm i0) \\
=& \chi_j +  [H_{free(j)},\chi_j]R(\lambda \pm i0).
\end{split}
\nonumber
\end{equation}

It follows from (\ref{12.09.2}), (\ref{12.09.6}) that 
\begin{equation}
\chi_jR(\lambda \pm i0) = R_{free(j)}(\lambda \pm i0)\chi_j  + 
R_{free(j)}(\lambda \pm i0)[H_{free(j)},\chi_j]R(\lambda \pm i0).
\label{S4chijresolvent}
\end{equation}

Note that on $\mathcal M_j$,
$f =  [H,\chi_j]u_j^{(-)}$,
and $[H,\chi_j]$ is a 1st order differential operator with coefficients which are compactly supported in $\mathcal M_j$.
Therefore, $f_j$ is compactly supported, and
\begin{equation}
\chi_jR(k^2 + i0)f = R_{free(j)}(k^2 + i0)f_j.
\label{S4chijandfj}
\end{equation}
Representing the resolvent in terms of modified Bessel functions 
(see (\ref{S5G0yy'})), if $m \neq 0$, we have for large $y$
\begin{equation}
\begin{split}
& \langle \chi_jR_{free(j)}(k^2 + i0)f_j,\varphi_{j,m}\rangle \\
&=  y^{(n-1)/2}K_{-ik}(\sqrt{\lambda_{j,m}}y)
\int_0^y(y')^{(n-1)/2}I_{-ik}({\sqrt{\lambda_{j,m}}}y')f_{j,m}(y')
\frac{dy'}{(y')^n}.
\end{split}
\label{S6Rfree(j)varphijm}
\end{equation}
Note that $K_{-ik}(z) = K_{ik}(z)$.
The case $m=0$ is computed similarly. We have thus proven the lemma. \qed

\medskip
Given $u_j^{(-)}$, $j=1,\cdots,N$, one can compute $b_{j,m}$ by observing the asymptotic behavior of 
$u - u^{(j)}$ in a neighborhood of the cusp. With this in mind, we make the following definition.


\begin{definition} \label{def:7.4}
We call the operator
$$
{\mathcal S}(k) : {\bf A}_{-\infty} \ni \psi^{(-)} \to \psi^{(+)} \in {\bf A}_{\infty}
$$
the {\it generalized S-matrix}.
\end{definition}

\begin{remark} 
 Definition 6.4 of the {\it generalized S-matrix} remains valid 
in the case that the metric on the cusp is not of the direct product form, but is perturbed by a super-exponentially decaying term. 
For example, we can assume that (\ref{S1ds2cuspexpand}), $1 \leq j \leq N$, is replaced by
(\ref{S1ds2regularexpand}) with coefficients $a_{j,pq}(\tilde x,y),b_{j,p}(\tilde x,y), c_{j}(\tilde x,y)$ satisfying the condition 
\begin{equation}
\Big((y\partial_{\tilde x})^{\alpha}(y\partial_y)^{\beta}\,a_{j,pq}(\tilde x,y)\Big)
e^{Cy} \in \mathcal B(M_j\times(1,\infty)).
\end{equation}
Here $C > 0$ is arbitrarily, and a similar estimate is valid for $b_{j,p}$, $c_j$. In fact, due to the formulae (\ref{S4chijresolvent}), (\ref{S6Rfree(j)varphijm}), the proof of Lemma \ref{S6LongLemma} remains valid for this case. 
\end{remark}


\subsection{Splitting the manifold}
We split $\mathcal M$ as
\begin{equation} \label{12.09.7}
\mathcal M = \mathcal M_{ext}\cup\mathcal M_{int}, \quad 
\mathcal M_{ext}= M_1 \times (2, \infty), \quad \mathcal M_{int} =\mathcal M \setminus 
\left(M_1 \times [2, \infty)  \right) ,
\end{equation}
Thus, $\mathcal M_{ext}$ and $\mathcal M_{int}$ have common boundary $\Gamma= M_1 \times \{2\}$.
Recall that the end $\mathcal M_1$ has a cusp, and $\mathcal M_{ext}$ is a direct
product of $M_1$ and $(2, \infty)$, i.e.
$$
(ds)^2= y^{-2} \left( (dy)^2+ h_1(x, dx) \right).
$$

Let $\Delta_g$ be the Laplace-Beltrami operator on $\mathcal M$, $H_{ext}$ and $H_{int}$ be 
$- \Delta_g - (n-1)^2/4$ defined on $\mathcal M_{ext}$, $\mathcal M_{int}$ with Neumann boundary 
condition on $\Gamma$, respectively. If $\mathcal M$ has only one end 
(i.e. $ N+N' = 1$), 
$\mathcal M_{int}$ is a compact manifold, and $H_{int}$ has a discrete spectrum. If $N \geq 2$, 
both of $\mathcal M_{int}$ and $\mathcal M_{ext}$ are non-compact, and the theorems in \S 4 
and \S 5 also hold in this case. We denote the inner product of $L^2(\Gamma)$ by
$$
\langle f,g\rangle_{\Gamma} = \int_{\Gamma}f\overline{g}\,d\Gamma.
$$

We put  
\begin{equation}
\Phi_m^{(0)} = \left\{
\begin{split}
& y^{(n-1)/2-ik} |M_1|^{-1/2}, \quad m = 0,\\
& y^{(n-1)/2}I_{-ik}({\sqrt{\lambda_{1,m}}}y) \varphi_{1,m }(x), \quad m \neq 0,
\end{split}
\right.
\nonumber
\end{equation}
\begin{equation}
g_m = (H-k^2)\chi_1\Phi_m^{(0)} = [H_{free(1)},\chi_1]\Phi_m^{(0)},
\nonumber
\end{equation}
\begin{equation}
\Phi_m =  \chi_1\Phi_m^{(0)} - R(k^2 + i0)g_m.
\label{S4phin}
\end{equation}


\begin{lemma} \label{L:7.5}
Let $k > 0$ and  $k^2 \not\in \sigma_p(H)\cap\sigma_p(H_{int})$. If $f \in L^2(\Gamma)$ satisfies
\begin{equation}
\langle f, \partial_{\nu}\Phi_{m}\rangle_{\Gamma} = 0, \quad \forall m \in {\bf Z}_+,
\label{S4forthdelyphi}
\end{equation}
where $\nu$ is the unit normal to $\Gamma$,
then $f = 0$.
\end{lemma}

Proof. 
We define an operator $\delta_{\Gamma}' \in {\bf B}(H^{-1/2}(\Gamma);H^{-2}(\mathcal M))$ by
$$
(\delta_{\Gamma}'v,w) = \langle v,\partial_{\nu}w\rangle_{\Gamma}, \quad 
\forall v \in H^{-1/2}(\Gamma), \quad 
\forall w \in H^2(\mathcal M),
$$
and define $u = R(k^2 - i0)\delta_{\Gamma}'f$ by duality, i.e. for $w \in L^{2,s}$, $s > 1/2$,
\begin{eqnarray*}
(R(k^2 - i0)\delta_{\Gamma}'f,w) &=& (\delta_{\Gamma}'f,R(k^2 + i0)w) \\
&=& \langle f,\partial_yR(k^2 + i0)w\rangle_{\Gamma}. 
\end{eqnarray*}
Then, if $f= \sum {\widehat f}_m \varphi_{1,m}$,
we have
\begin{equation}
R(k^2 -i0)\delta_{\Gamma}'f = 
\sum_{m\in{\bf Z}_+}A_m(y)\widehat f_m \varphi_{1,m},
\label{S4R0Ndeltagamma}
\end{equation}
where for $m \geq 1$,
\begin{equation}
A_m(y) = \left\{
\begin{split}
\left(y^{(n-1)/2}K_{ik}({\sqrt{\lambda_{1,m}}} y)\right)'
\Big|_{y=2}y^{(n-1)/2}I_{ik}(\sqrt{\lambda_{1,m}}y), 
\quad 
y < 2, \\
\left(y^{(n-1)/2}I_{ik}({\sqrt{\lambda_{1,m}}}y)\right)'\Big|_{y=2}y^{(n-1)/2}
K_{ik}({\sqrt{\lambda_{1,m}}} y), 
\quad 
y > 2, 
\end{split}
\right.
\label{S4An}
\end{equation}
and for $m = 0$,
\begin{equation}
A_0(y) = \left\{
\begin{split}
\left(y^{(n-1)/2-ik}\right)' 2^{(n-1)/2+ik}, 
\quad 
y < 2, \\
\left(y^{(n-1)/2+ik}\right)' 2^{(n-1)/2-ik}, 
\quad 
y > 2.
\end{split}
\right.
\label{S4A0}
\end{equation}
Then 
$(H - k^2)u = \delta_{\Gamma}'f$ in the sense of distribution, hence
\begin{equation}
(H- k^2)u = 0 \  {\rm except \ for} \ \Gamma.
\label{S4zerooutsidegammma}
\end{equation}
By  (\ref{S5G0yy'}) and (\ref{S4chijresolvent}), we have when $y$ is large enough and $m \neq 0$,
\begin{eqnarray*}
\langle \chi_Nu,\varphi_{1,m}\rangle &=& 
y^{(n-1)/2}K_{ik}({\sqrt{\lambda_{1,m}}}y)
\int_0^{\infty}\int_{M_1}\varphi_{1,m}(x)  (y')^{(n-1)/2}I_{ik}({\sqrt{\lambda_{1,m}}}y')\\
& &\ \ \ \times \left\{\chi_1 + [H_{free(1)},\chi_1]R(k^2 + i0)\right\}
\delta_{\Gamma}'f\frac{dy'}{(y')^n} \\
&=& y^{(n-1)/2}K_{ik}({\sqrt{\lambda_{1,m}}}y)\left(\delta_{\Gamma}'f,\left\{\chi_1 - R(k^2 + i0)[H_{free(1)},\chi_1]\right\}\Phi_m^{(0)}\right) \\
&=& y^{(n-1)/2}K_{ik}({\sqrt{\lambda_{1,m}}}y)\langle f,\partial_y\Phi_m\rangle_{\Gamma} = 0.
\end{eqnarray*}
Similarly, one can show for large $y$ that
$$
\langle\chi_1u,\varphi_{1,0}\rangle = 0.
$$
Therefore $u = 0$ when $y$ is large enough. Since $(H - k^2)u = 0$ in $\mathcal M_{ext}$, the unique continuation theorem imply that $u = 0$ in $\mathcal M_{ext}$. By the resolvent equation (\ref{S4chijresolvent}), formulae 
(\ref{S4R0Ndeltagamma}) $\sim$ (\ref{S4A0}), and the fact 
that $[H_{free(1)}, \chi_1] R(\lambda+i0) \delta'f$  is smooth,
we see that $\partial_yR(k^2 - i0)\delta_{\Gamma}'f$ is continuous across $\Gamma$. Therefore, in $\mathcal M_{int}$, $u$ satisfies $(H_{int} - k^2)u = 0$ and the Neumann boundary condition on $\Gamma$, hence $u = 0$ in $\mathcal M_{int}$. This follows from the assumption $k^2 \not\in \sigma_p(H_{int})$ when $\mathcal M_{int}$ is  compact, and from Lemma \ref{S4RadCondUniqueness} when $\mathcal M_{int}$ is non-compact. We thus have $u = 0$ in $\mathcal M$, and $f = 0$. \qed

\medskip
The generalized S-matrix ${\mathcal S}(k)$ is an operator-valued 
$(N+N')\times (N+N')$ matrix. Let ${\mathcal S}_{11}(k)$ be its $(1,1)$ entry. 
For $a  \in \ell^{2,\infty}$, we put
$b = {\mathcal S}_{11}(k)a \in \ell^{2,-\infty}$, and
$$
\Phi = \sum_{m\in{\bf Z}_+}a_m\Phi_m.
$$
Then $(H-k^2)\Phi = 0$ and in $\mathcal M_1$, it takes the form
$$
\Phi = u_1^{(-)} - u_1^{(+)},
$$
$$
u_1^{(-)} = \omega_-^{(c)}(k)\Big(a_{0}y^{(n-1)/2-ik} + 
\sum_{m\geq 1}a_m \varphi_{1,m}(x) y^{(n-1)/2}I_{-ik}
({\sqrt{\lambda_{1,m}}}y)\Big),
$$
$$
u_1^{(+)} = \omega_+^{(c)}(k)\Big(b_{0}y^{(n-1)/2+ik} + 
\sum_{m\neq 0}b_m \varphi_{1,m}(x) y^{(n-1)/2}K_{ik}({\sqrt{\lambda_{1,m}}}y)\Big).
$$
Therefore, to determine ${\mathcal S}_{11}(k)$ is equivalent 
to the observation of the outgoing  exponentially decaying waves $u_1^{(+)}$ at $\mathcal M_1$
for any
 incoming exponentially growing waves $u_1^{(-)}$ at the cusp $\mathcal M_1$. 


\subsection{Gel'fand problem, BSP and N-D map}
Before going to proceed, let us recall the {\it Gel'fand problem}. Let $\Omega$ be a compact 
Riemannian manifold with boundary $\Gamma = \partial\Omega$, and $- \Delta_g$ the associated Laplace-Beltrami operator. Let $0=\lambda_1 < \lambda_2 < \cdots$ be its Neumann eigenvalues without counting multiplicities, and $\varphi_{i,1},\cdots,\varphi_{i,m(i)}$ be the orthonormal system of eigenvectors associated with the eigenvalue $\lambda_i$. Let us call the set
$$
\Big\{(\lambda_i,\varphi_{i,1}\big|_{\Gamma},\cdots,\varphi_{i,m(i)}\big|_{\Gamma})\Big\}_{i=1}^{\infty}
$$
the boundary spectral data ({\bf BSD}).
The problem raised by Gel'fand is :

\begin{quote}
{\it Does BSD determine the Riemannian metric of $M$?}
\end{quote}

This problem was solved by Belishev-Kurylev \cite{BeKu92} using the boundary control method (BC-method) proposed by 
Belishev \cite{Be87}. 
The BC method was advanced in \cite{KKL01} and we use in this paper this variant of the method. On different variations if
this technique, see  \cite
{AKKLT,BKLS,KKL08}.

Although it is formulated in terms of BSD, what is actually used in 
the BC-method is the boundary spectral projection ({\bf BSP}) defined by
\begin{equation}
\Big\{(\lambda_i, \sum_{j=1}^{m(i)}\varphi_{i,j}(x)\overline{\varphi_{i,j}(y)}\big|_{(x,y)\in\Gamma\times\Gamma})\Big\}_{i=1}^{\infty}.
\label{S4BSPdiscrete}
\end{equation}
This appears in the kernel of the Neumann to Dirichlet map (N-D map)
\begin{equation}
\Lambda(z) : f \to u,
\label{S4NDmap}
\end{equation}
where $u$ is the solution to the Neumann problem
\begin{equation}
\left\{
\begin{split}
& (- \Delta_g - z)u = 0 \quad {\rm in} \quad \Omega, \\
& \partial_{\nu}u = f \in H^{-1/2}(\Gamma),
\end{split}
\right.
\label{S4Neumannproblem}
\end{equation}
$\nu$ being the outer unit normal to $\Gamma$, $z \not\in \sigma(-\Delta_g)$.
The N-D map is related to the resolvent $(- \Delta_g - z)^{-1}$ in the following way :
\begin{equation}
\Lambda(z) = \delta_{\Gamma}^{\ast}(- \Delta_g - z)^{-1}\delta_{\Gamma}, \quad 
z \not\in \sigma(-\Delta_g),
\label{S4NDreslvent}
\end{equation}
where $\delta_{\Gamma} \in {\bf B}(H^{-1/2}(\Gamma);H^{-1}(\Omega))$ is the adjoint of the trace operator 
$$
r_{\Gamma} : H^1(\Omega) \ni w \to w\big|_{\Gamma} \in H^{1/2}(\Gamma),
$$
\begin{equation}
(\delta_{\Gamma}f,w)_{L^2(\Omega)} = (f,r_{\Gamma}w)_{L^2(\Gamma)}, 
\quad f \in H^{-1/2}(\Gamma), \quad w \in H^1(\Omega).
\label{S4deltagamma}
\end{equation}
In our case, we have

\begin{lemma} \label{L:7.6}
In the case of $\Omega= {\mathcal M}_{int}$ of form (\ref{12.09.7}),
to give BSP is equivalent to give the N-D map $\Lambda(z)$ for all $z \not\in \sigma(-\Delta_g)$.
\end{lemma}

We need to exlain more about BSP. If $\overline{{\mathcal M}_{int}}$ is a compact manifold with boundary, the BSP defined abov works well.
If $\Omega={\mathcal M}_{int}$ is non-compact with compact boundary $\Gamma$, then its shifted
Laplace operator $H_{int} = -\Delta_g-(n-1)^2/4$ with Neumann boundary condition has a continuous spectrum 
$\sigma_c(H_{int}) = [0,\infty)$, and, furthermore, $H$ has a spectral representation $\mathcal F$ 
like the one discussed in \S5. Then the notion of BSP should be modified as 
\begin{equation}
\Big\{\delta_{\Gamma}^{\ast}\mathcal F(k)^{\ast}\mathcal F(k)\delta_{\Gamma}\, ; \, k > 0\Big\}\cup\Big\{(\lambda_i,\delta_{\Gamma}^{\ast}P_i\delta_{\Gamma})\Big\}_{i=1}^{m},
\nonumber
\end{equation}
where $\lambda_i$ is the eigenvalue of $H_{int}$, $P_i$ is the associated eigenprojection 
and $m$ is the number, finite or infinite, of eigenvalues.
In this case, we extend the N-D map $\Lambda(z)$ for $z = k^2 \in (0,\infty)\setminus\sigma_p(H)$ 
by using the solution $u$ of (\ref{S4Neumannproblem}) satisfying the outgoing radiation condition. 
Then Lemma \ref{L:7.6}  holds true. 
(See Lemma 3.3 in \S5, \cite{IsKu09} and Lemma 5.6 in \S 4, \cite{IKL10}.)

For another approach to the manifold reconstruction, which uses N-D map at one frequency only, see e.g. \cite{LU01} and \cite{LeU89}.


\subsection{Generalized S-matrix and N-D map}

Returning to our problem, we define the N-D map for $\mathcal M_{int}$ by (\ref{S4NDmap}) 
and (\ref{S4Neumannproblem}). 

Now suppose we are given two orbifolds $\mathcal M^{(i)}$, $i = 1,2$, satisfying the 
assumptions (A-1) $\sim$ (A-3) in \S 1. Let $H^{(i)}$ be the Laplace-Betrami operator of 
$\mathcal M^{(i)}$. Assume that $\mathcal M^{(i)}$ has $N_i+N_i'$ numbers of ends,
$N_i>0$, and let 
${\mathcal S}^{(i)}(k)$ be the $(1,1)$ entry of the generalized S-matrix for $H^{(i)}$. Let 
$({\mathcal M}^{(i)}_{1}, h_1^i), i=1,2,$ be isometric Riemannian orbifolds. Thus, we can naturally 
identify ${\mathcal M}_{1}^{(i)}, i=1,2,$. Next we split 
$\mathcal M^{(i)}$ into $\mathcal M_{int}^{(i)} \cup \mathcal M_{ext}^{(i)}$ as above using 
$\Gamma^{(1)}=\Gamma^{(2)}$. Let $H_{int}^{(i)}$ be the Laplace-Beltrami operator of 
$\mathcal M_{int}^{(i)}$ with the
Neumann boundary condition on $\Gamma^{(i)}$, and define the N-D map $\Lambda^{(i)}(z)$ for 
$\mathcal M_{int}^{(i)}$. With this preparation, we can prove the following lemma.

\begin{lemma} \label{L:7.7}
If ${\mathcal S}^{(1)}_{11}(k) = {\mathcal S}^{(2)}_{11}(k)$ for  $k > 0$, $k^2 \not\in \sigma_p(H^{(1)})\cup\sigma_p(H^{(2)})$, we have 
$$
\Lambda^{(1)}(k^2) = \Lambda^{(2)}(k^2), \quad \forall k^2 \in (0,\infty)\setminus\sigma_p(H^{(1)}_{int})\cup\sigma_p(H^{(2)}_{int}).
$$
 Moreover, BSP's for $H_{int}^{(1)}$ and $H_{int}^{(2)}$ coincide.
\end{lemma}

Proof. For $i = 1, 2$, we construct $\Phi_n^{(i)}$ as in (\ref{S4phin}), and put 
$$
u = (\Phi_m^{(1)} - \Phi_m^{(2)})|_{{\mathcal M}_i}. 
$$
Then $u$ satisfies $(H^{(i)} - k^2)u = 0$ in 
$\mathcal M_{ext}^{(1)} = \mathcal M_{ext}^{(2)}$. Since ${\mathcal S}^{(1)}_{11}(k) = {\mathcal S}^{(2)}_{11}(k)$, 
by Definition \ref{def:7.4},  $u = 0$ in $\mathcal M_{ext}^{(1)} = \mathcal M_{ext}^{(2)}$. 
Hence, $\partial_{\nu}\Phi_m^{(1)} = \partial_{\nu}\Phi_m^{(2)}$ on $\Gamma$. 

In $\mathcal M_{int}^{(i)}$, $\Phi_n^{(i)}$ is the outgoing solution of the equation $(H^{(i)} - k^2)v = 0$. Hence, $ \Phi_n^{(i)}\big|_{\Gamma} = \Lambda^{(i)}(k)(\partial_{\nu} \Phi_n^{(i)})\big|_{\Gamma}$. 
This implies that
\begin{equation}
\Lambda^{(1)}(k)(\partial_{\nu}\Phi_m^{(1)})\big|_{\Gamma} = 
\Lambda^{(2)}(k)(\partial_{\nu}\Phi_m^{(2)})\big|_{\Gamma}, \quad \forall m.
\label{S4NDcoincide}
\end{equation} 
Lemma \ref{L:7.5} implies that the linear span of 
$\{\partial_\nu \Phi_m^{(i)}\big|_{\Gamma}\, ; \, m \in {\bf Z}\}_+$ is dense in $L^2({\Gamma})$. 
Therefore, by (\ref{S4NDcoincide}), $\Lambda^{(1)}(k) = \Lambda^{(2)}(k)$,
for $k>0,\, k \notin \cup_{i=1}^2 \sigma_p(H^{(i)}) \cup \sigma_p(H^{(i)}_{int})$.  \qed


\section{Orbifold isomorphism}

Above we have reduced the inverse scattering problem  for the construction of an orbifold from local measurements. 
This problem is studied in detail in  \cite{KLY09}  in the general $n$-dimensional case and in
 the 2-dimensional case, in the context of scattering problems, in \cite {IKL10}. 
We give an outline of the proof of Theorem 1.1. that is based on the basic steps of \S5 in \cite{IKL10}, and 
with modifications necessary to deal with the multidimensional, rather than 2D, case.

Observe thar $\Lambda(z)$ is an integral operator with the kernel $G_{int}(X, Y; z),\, x, x' \in \Gamma$,
where $G_{int}(\dot, \cdot; z)$ is  Green's function for for the Neumann problem in $\mathcal M_{int}$.
The function $G_{int}(X, Y; z)$ enjoys the following separation property:

\noindent If $G_{int}(X, Y; z)=G_{int}(X', Y; z)$ for all $Y \in \Gamma,\, z \in {\bf C}\setminus {\bf R}$,
then $X=X'$.

The proof of Theorem 1.1 consists of two principal steps. First, we show that 
${\mathcal M}_{int}^{(i), reg}$ are isometric and, therefore
${\mathcal M}_{int}^{(i)}$ are isometric as metric spaces. In the future, we refer to this isometry
as ${\bf X}:{\mathcal M}_{int}^{(1)} \to   {\mathcal M}_{int}^{(2)}$ and, in the first step, we recover
$\bf X$ in an inductive procedure. 

Next, we use the fact that, if two oriented orbifolds are isometric as metric spaces, then
they are orbifold isometric, see \S1.4. This effectively means that if 
$p^{(i)} \in {\mathcal M}^{(i), sing}$ and $p^{(2)}= {\bf X}(p^{(1)})$, then the groups $G_{p^{(1)}}$ and
$G_{p^{(2)}}$  are isomorphic.

The proof of the first step is made by inductively enlarging the parts of ${\mathcal M}_{int}^{(i), reg}$
which are isometric to each other by starting from $\Gamma=\Gamma^{(1)}=\Gamma^{(2)}$.
Let us call the parts of ${\mathcal M}_{int}^{(i), reg}$ which are already proven to be isometric after
$m$ iterations by $\Omega_m^{(i)}$. It would always be the case that $\Omega_m^{(i)}$ are arcwise connected
with $\Gamma$. Then, from Lemma \ref{L:7.7}, we see that
$$
G_{int}^{(2)}({\bf X}(X), {\bf X}(Y); z)=G_{int}^{(1)}(X, Y; z),\quad X, Y \in \Omega_m^{(1)},
\quad z \in {\bf C} \setminus {\bf R}.
$$
The induction step is based on the possibility to extend the isometry ${\bf X}:\Omega_m^{(1)}
\to \Omega_m^{(2)}$ from 
$B_\delta(p^{(1)})= {\bf X}^{-1}(B_\delta(p^{(2)}),\, B_\delta(p^{(i)}) \subset \Omega_m^{(i)}$
to ${\bf X}: B_r(p^{(1)}) \mapsto B_r(p^{(2)})$ for any
$r\leq \min\{\hbox{inj}(p^{(1)}),\, \hbox{inj}(p^{(2)})\}$, where $\hbox{inj}(p)$ is the injectivity
radius at $p$. Since the geodesics stop when they hit a singular point, 
$\hbox{inj}(p^{(i)}) \leq d_i(p^{(i)},\, {\mathcal M}_{int}^{(i), sing})$.
Using arguments, based on the determination of the volume of balls
$B_\rho(x)\subset {\mathcal M}^{(i)}$, $x\in   B_r(p^{(i)})$, that is analogous to the proof of Theorem 5.7 
in \cite{IKL10}, or alternatively, using techniques developed in \cite{KirK},
we see that for  any $r\leq \min\{\hbox{inj}(p^{(1)}),\, \hbox{inj}(p^{(2)})\}$ 
the set  $\partial B_r(p^{(1)})\cap  {\mathcal M}_{int}^{(1), sing}$ is non-empty
if and only if 
 $\partial B_r(p^{(2)})\cap  {\mathcal M}_{int}^{(2), sing}$ is non-empty.
These imply, in particular, that $\hbox{inj}(p^{(1)})=\hbox{inj}(p^{(2)})$.
Thus, we can extend the isometry ${\bf X}:\Omega_m^{(1)}
\to \Omega_m^{(2)}$ with the map
 ${\bf X}: B_{\rm {inj}(p^{(1)})}(p^{(1)}) \mapsto B_{\rm{inj}(p^{(2)})}(p^{(2)})$. 
%
%
%
%
%
%
%
%
This implies we can avoid  the set
${\mathcal M}_{int}^{(i), sing}$ when we extend the isometry ${\bf X}$
to a larger set.
Hence,
using the fact that ${\mathcal M}_{int}^{(i), reg}$ are path-connected,
we can form an iterative procedure where the set $\Omega_m^{(i)}\subset 
{\mathcal M}_{int}^{(i), reg}$ is made larger in each step,. By
constructing a maximal such extension we can extend $\bf X$ from $\Gamma$ onto ${\mathcal M}_{int}^{(1), reg}$ .

The second step, regarding the orbifold isomorphism of ${\mathcal M}_{int}^{(i)}$ can be carried
out following the same considerations as in the proof of Theorem 9.1 in \cite{KLY09}. Note
that, although Theorem 9.1 there is proven only for the compact orbifolds, it can be easily
extended to the case of the non-compact orientable orbifolds. Indeed, the constructions
in the proof of Theorem 9.1 are local dealing only with the neighbourhoods of the singular points
in the orbifold.
\qed

\end{document}